\documentclass{article}
\usepackage[utf8]{inputenc}
\usepackage{amsfonts}
\usepackage{lscape}
\usepackage{amsmath}
\usepackage{amssymb}
\usepackage{braket,amsfonts,amsopn}
\usepackage{amsfonts}

\usepackage{amsthm}
\usepackage{todonotes}
\usepackage{hyperref}

\newcommand{\R}{\mathbb R}
\newcommand{\N}{\mathbb N}
\newcommand{\E}{\mathbb E}
\newcommand{\eps}{\varepsilon}

\newcommand{\1}[1]{\chi_{#1}}
\newcommand*{\defeq  }{\mathrel{\vcenter{\baselineskip0.5ex\lineskiplimit0pt\hbox{\normalsize.}\hbox{\normalsize.}}}=}
\newcommand{\Frechet}{Fr\'echet}
\newtheorem{lemma}{Lemma}

\theoremstyle{plain}
\newtheorem{rem}{Remark}
\newtheorem{algo}{Algorithm}

\newcommand{\review}[1]{{{#1}}}

\title{Wavelet-based priors accelerate maximum-a-posteriori optimization in Bayesian inverse problems} 
\author{Philipp Wacker\and Peter Knabner\\Friedrich-Alexander-Universit\"at Erlangen-N\"urnberg, }

\begin{document}
\maketitle
\begin{abstract}
Wavelet (Besov) priors are a promising way of reconstructing indirectly measured fields in a regularized manner. We demonstrate how wavelets can be used as a localized basis for reconstructing permeability fields with sharp interfaces from noisy pointwise pressure field measurements \review{in the context of the elliptic inverse problem}. For this we derive the adjoint method of minimizing the Besov-norm-regularized misfit functional (this corresponds to determining the maximum a posteriori point in the Bayesian point of view) in the Haar wavelet setting. As it turns out, choosing a wavelet--based prior allows for accelerated optimization compared to established trigonometrically--based priors. 
\end{abstract}
\section{Introduction}
We imagine we want to infer some unknown (typically infinite-dimensional) quantity $u$ by their effect (through some ``forward map'') on a noisy observed quantity $y$ (in general finite-dimensional), i.e.
\begin{equation}\label{eq:invProb}
y = G(u) + \eta
\end{equation}
where $\eta$ is observational noise. By dimensionality and randomness, this is a vastly ill-posed problem in dire need of some sort of regularization. For a classical regularization approach, see \cite{engl1996regularization}, or in the context of groundwater inference, \cite{sun1985identification}. In the Bayesian context, this means choosing a prior model describing our belief concerning a-priori feasible values of $u$. Then, Bayes' theorem allows us to incorporate the observation $y$ in order to obtain the posterior distribution $u|y$. 
In many applications, the inverse problem consists of inferring a function (i.e. a permeability or density scalar field over some domain) by their effect on some observation through a PDE-governed dynamics. This raises the question, ``what are random functions?'' 

This is a very subtle question and even its simplification, ``what are Gaussian random functions?'', is very difficult to answer (see \cite{bogachev1998gaussian} and \cite{kuo1975gaussian}) and full of pitfalls due to the infinite dimensionality of function spaces. 

In the Gaussian setting, there are basically two ways of defining random functions: Either by interpreting them as samples from a Gaussian random field which can be uniquely defined by its mean field and covariance function, or by defining an abstract Gaussian measure on a vector space of functions, uniquely determined by a mean function and a covariance operator. The latter view intuitively corresponds to defining a random function as a randomized superposition of deterministic basis functions of the function space (e.g. a trigonometric basis, this corresponds to the Gaussian prior in section \ref{sec:commonprior}). This is due to the fact that any Gaussian measure on a Banach space has a white noise expansion (see \cite[Theorem 2.12]{da2014stochastic}).

A relatively new (non-Gaussian) function space prior is the Besov prior (which we will call Wavelet prior in this paper, for reasons stated below) described mathematically in \cite{dashti2011besov}, \cite{lassas2009discretization} and practically demonstrated in one dimension in \cite{bui2015scalable,kolehmainen2012sparsity}. The main idea of this prior is that a random function sampled from this prior is given by a random superposition of localized basis functions (a set of compactly supported wavelets) in constrast to the non-local trigonometric basis functions often employed.

This article is supposed to showcase the application of Besov priors in Bayesian inverse problems for higher dimensions, in particular $d=2$. The main contribution of this article is the following: When computing the maximum-a-posteriori point (MAP) in its (already more efficiently solvable) adjoint formulation, it turns out that using a wavelet--based prior facilitates an important computational speed-up in a critical step of the optimization procedure. This is explained in more detail in remark \ref{rem:main},

In a few numerical experiments we show how Wavelet priors can be employed to reconstruct essential information about permeability fields with sharp interfaces in the inverse elliptical problem setting. This is in contrast to the established $\beta\cdot(-\Delta)^{\alpha/2}$--prior, which has essentially two drawbacks in our setting: Reconstructions using Gaussian priors based on trigonometric expansions will invariably lead to smoothened reconstructions of jumps. Also, the computational speed-up present with Wavelet-based inversion is not available.
In this study, we will restrict ourselves to the Haar wavelet basis so that we can use its discontinuity to model interfaces.

Wavelets have been employed in a groundwater inference setting in \cite{chatterjee2012multi}, but their approach differs in philosophy from the Bayesian wavelet prior notion in that they employ a template matching algorithm using training images. In \cite{wang2017bayesian}, a fully Bayesian approach with sampling from the posterior was presented in the framework of the Randomize-then-optimize method, but we focus here on computing the maximum-a-posteriori estimator. The reason for this is that the wavelet setting allows for a very efficient way of computing this quantity using the adjoint method. 

There is also recent progress in the field of medical imaging using Besov priors \cite{rantala2006wavelet}. The use of wavelets in inverse problems is not entirely new, see for example \cite{rieder1997wavelet,dubot2015wavelet}, but here wavelets were used in the context of multilevel methods. In \cite{wang2017bayesian}, wavelet-based Besov priors were employed in the demonstration of the RTO method of sampling.

The rest of the article is structured as follows: The rest of this section gives a short introduction into (Bayesian) inverse problems, section 2 treats the model inverse problem we use and describes a common prior used in practice and how it is unfavorable here. Section 3 gives a short introduction into the construction of Besov prior using Haar wavelets. The main section 4 presents how wavelets and the adjoint approach in optimization go hand in hand to allow for efficient computation of the maximum-a-posteriori point. Section 5 presents a few numerical simulations for illustration and section 6 concludes with a few open directions and ideas for improvements.

\subsection{General Bayesian inverse problems}
The inverse problem we are interested in is defined as follows\footnote{we follow the exposition of the inverse problem in \cite{blscwa}}.
The goal is to recover the unknown parameters $u \in X$ from a noisy observation $y \in Y$ in \eqref{eq:invProb} with the so-called forward response map $G:X\to Y$ which often -- and in our case -- is a PDE solution operator. We denote by $X$ and $Y$ separable Hilbert spaces and by the random variable $\eta$ the observational noise.

Inverse Problems arise in a multitude of applications, for example in geosciences, medical imaging, reservoir modelling, astronomy and signal processing. A probabilistic approach has been undertaken first by Franklin \cite{franklin1970well}, Mandelbaum \cite{mandelbaum1984linear} and others, and a formulation of inverse problems in the context of Bayesian probability theory was originally obtained by Fitzpatrick \cite{fitzpatrick1991bayesian}. More recent literature about inverse problems in the Bayesian setting can be found in \cite{neubauer2008convergence}, \cite{stuart2010inverse}, \cite{helin2015maximum}. See also \cite{carrera1986estimation} for probabilistic groundwater estimation from a more applied viewpoint.

In the following, due to the nature of the inverse problem we are interested in, we assume that the number of observations is finite, i.e. $Y = \R^K$ for some $K\in\N$, whereas the  the unknown is a distributed quantity (and hence infinite--dimensional), which is a typical setting for many applications mentioned above. Thus, the inverse problem consists of inferring infinitely many parameters from a finite number of noisy observations. 
This leads to an ill-posed problem (this is also often true for the case of infinitely many observations, as in the case of the inverse heat equation, see \cite{stuart2010inverse}), which needs to be regularized. We will focus here on the Bayesian approach, i.e. the inverse problem becomes well-posed by a prior distribution on the feasible values for $u$ and the data gets incorporated by a Bayesian update on this prior measure, 
yielding the posterior measure in result.  Furthermore, we assume that the observational noise is Gaussian, i.e. $\eta\sim N(0,\Gamma)$ with a symmetric, positive definite matrix $\Gamma\in\mathbb R^{K\times K}$.

We consider the least squares ``error'' (or model-data misfit) functional
\begin{equation}
\label{eq:errorfunctional}
\Phi(u; y) = \frac{1}{2}\| \Gamma^{-1/2}(y-G(u))\|_Y^2,
\end{equation}
where $\Gamma$ normalizes the model-data misfit and is normally chosen as the covariance operator of the noise.
Plain infimization of this cost functional is not feasible due to the ill-posedness of the problem. For a given prior $u \sim \mu_0$, we derive the posterior measure $u|y \sim \mu$ where 
\begin{equation}\label{eq:posteriormeasure}
\mu(du) = \frac{1}{Z}\exp(-\Phi(u;y))\mu_0(du).
\end{equation}
For a derivation of this fact, see \cite{dashti2013bayesian}.
An issue of Bayesian inverse problems is how to obtain useful information from the posterior, which is in general given by \eqref{eq:posteriormeasure}, i.e. abstractly as a density w.r.t. the prior. Most quantities of interest can be computed by sampling, for example with Metropolis-Hastings methods. 

The \textit{maximum-a-posteriori point } plays a key role as a ``point of maximum density'' (although this does not make rigorous sense in the infinite dimensional setting). It is defined by the position of small balls asymptotically maximizing  posterior probability. In some cases (as for the Besov prior under consideration here) the MAP estimator can be shown  to be equivalently characterized as a minimizer of a ceratin functional. In the case where the prior measure is Gaussian with Cameron-Martin space $(E,\|\cdot\|_E)$, this is
\[u_{\text{MAP}} \defeq \operatorname{arginf} I(u) \]
where $I(u) = \Phi(u) + \frac{1}{2}\|u\|_E$. 

For a discussion and development of this subtle issue, see \cite{dashti2013map,helin2015maximum,agapiou2017sparsity}

Our simulations will mostly show the computed value of $u_{\text{MAP}}$. Another interesting quantity is the \textit{conditional mean}, which is
\[u_{\text{CM}} \defeq \E^\mu u.\]
In general, $u_{\text{CM}}$ and $u_{\text{MAP}}$ can be vastly different (but they are related from a theoretical perspective, see \cite{burger2014maximum,burger2016bregman}). Typically, the MAP point can exhibit sparse behaviour (if the corresponding optimization problem warrants this), but the Conditional Mean (as well as samples of the posterior) will often fail to have this property. This serves as a general caveat to ``sparse Besov priors''.\footnote{thanks to Youssef Marzouk who pointed this out to me.}

Note that we do not experiment with full posterior sampling here as our focus is the computational speed-up we experience with MAP optimization. In principle, though, this could be used for posterior sampling MCMC methods which use likelihood gradient information.
\section{The model inverse elliptical problem}\label{sec:model}
Our model inverse problem will be the following: Given a solution $p$ of 
\begin{align*}
-\nabla\cdot (e^{u(x, y)} \cdot \nabla p(x, y)) &= f(x, y)\text{ in $[0,1]^2$}\\
p(x, y) &= g(x, y) \text{ on $\Gamma_D$}\\
\partial_\nu p(x, y) &= h(x, y) \text{ on $\Gamma_N$}
\end{align*}
($\Gamma_D\cup\Gamma_N = \partial [0,1]^2$ is the decomposition of the boundary in Dirichlet and Neumann boundary) to some permeability field $\theta = e^u \in L^\infty([0,1]^2)$ (the field $u$ is called the logpermeability), we assume that we have noisy pointwise observations $y_i = p(x_i,y_i) + \eta_i$, $i=1,\ldots,N$ with $\eta_i\sim N(0,\gamma^2)$ i.i.d. The inverse problem is the inference
\[ \R^N \ni (y_i)_{i=1}^N \mapsto u\in L^\infty([0,1]^2),\]
i.e. we want to find the unknown logpermeability field $u$ from pointwise observations of the corresponding pressure field. The inference of an infinite-dimensional quantity from finite-dimensional one is obviously ill-posed, hence the need for a regularization, in our case a prior.

In the language of the last section, $\Phi(u; y) = \frac{1}{2\gamma^2}\|y - \Pi p\|^2$ where \[\Pi p = (p(x_1,y_1),\ldots,p(x_N, y_N))^T\] is the evaluation operator and $y = (y_i)_{i=1}^N$ is the observation vector. The survey paper by Delhomme et al \cite{delhomme2000four} and Sun's book about inverse problems \cite{sun2013inverse} give a thorough introduction to the evolution of inverse problems in the geohydrology community. 
\subsection{The trigonometric prior \texorpdfstring{$N(\mu, \beta(-\Delta)^{-\alpha})$}{N(mu,beta*(-Laplace)**(-alpha)}}\label{sec:commonprior}
\begin{figure}[hbtp]\label{fig:trigvswavelet}
\centering
\includegraphics[width=0.47\textwidth]{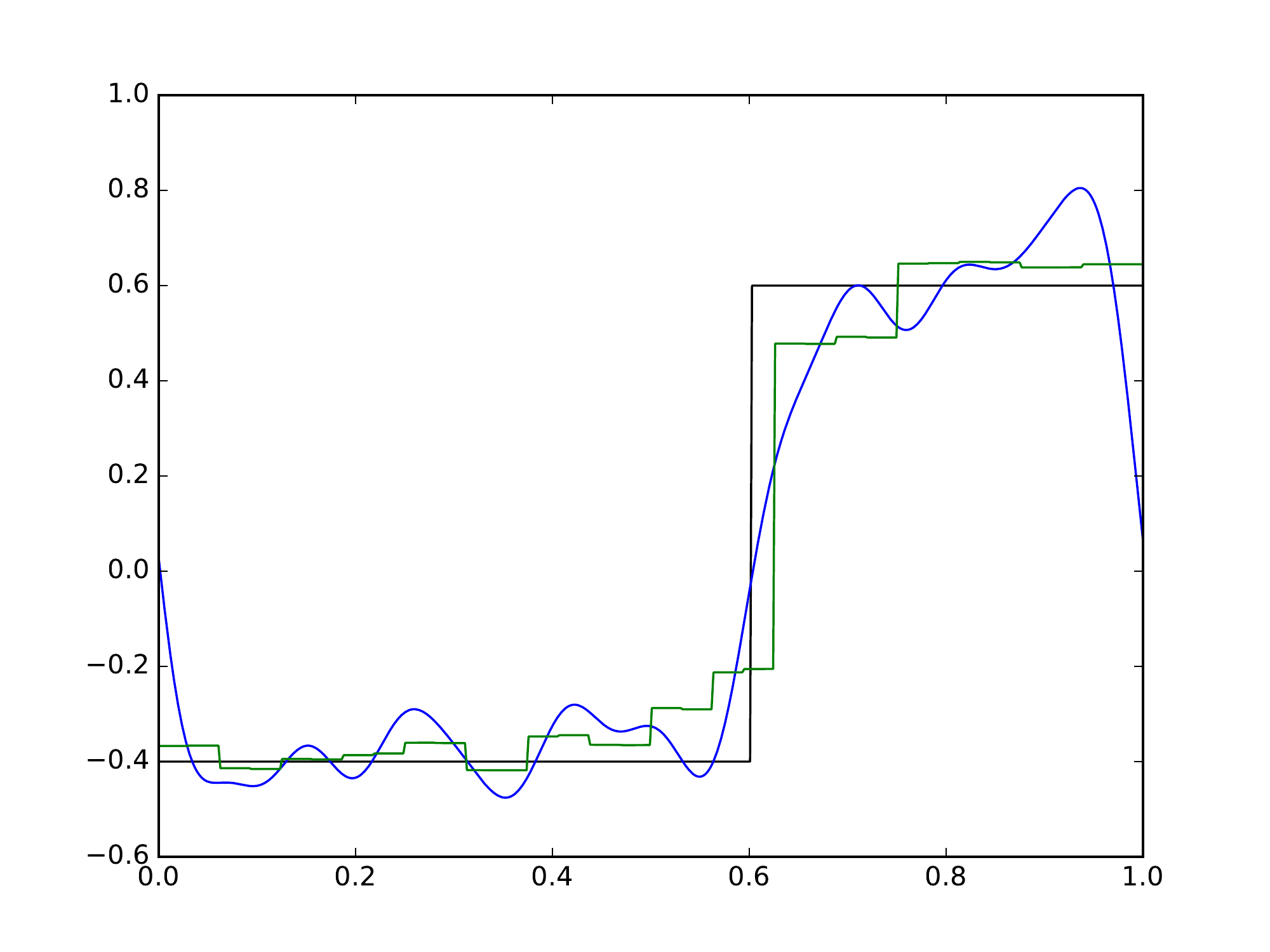}
\includegraphics[width=0.47\textwidth]{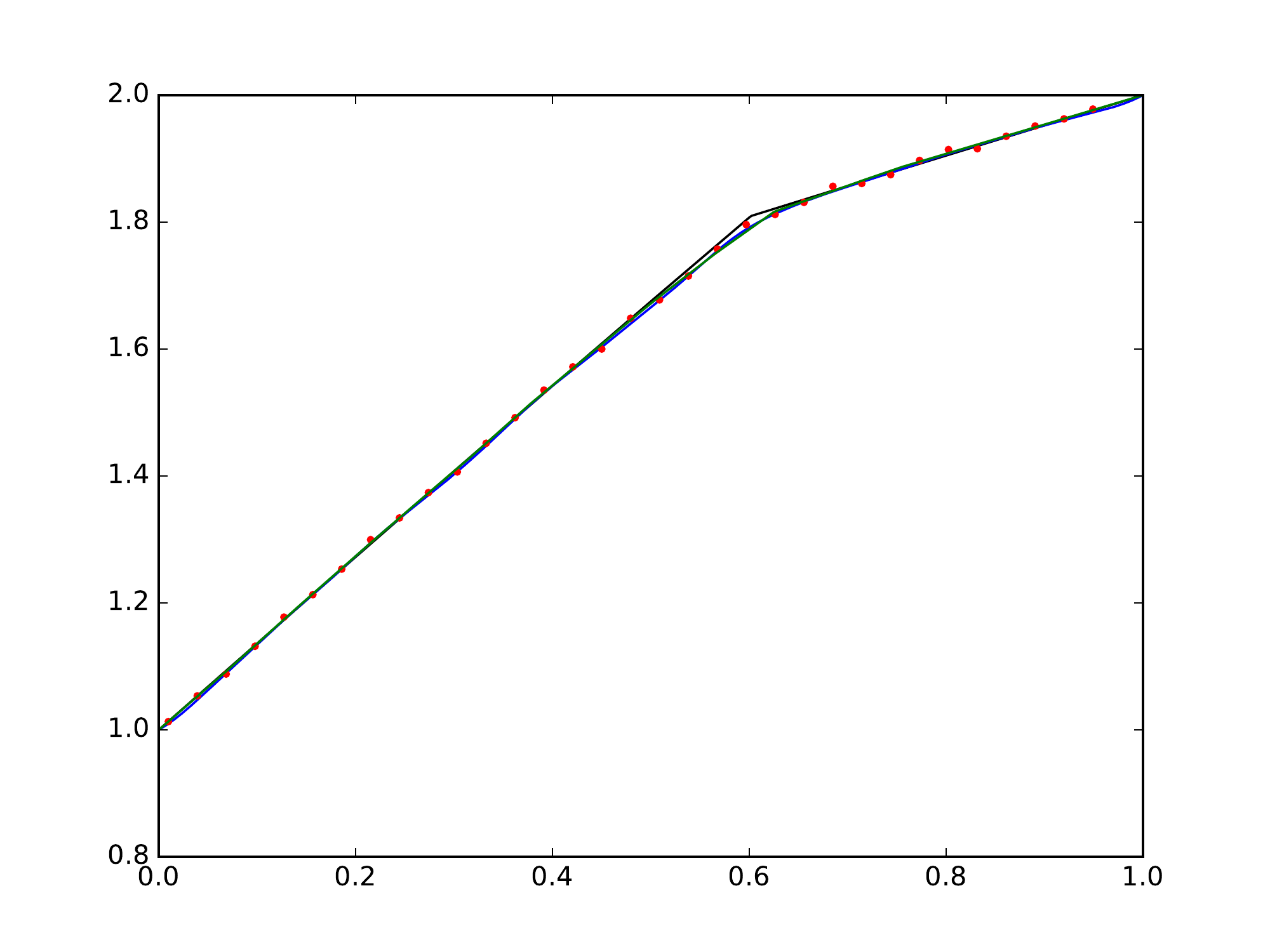}
\caption{Comparison of trigonometric and wavelet prior for a 1d inverse elliptical problem. Right: Pressure with observations. Left: Ground truth permeability (black piecewise constant with values $-0.4$ and $0.6$), $\mathrm{MAP}$ point for $(-\Delta)^{-1}$ prior (blue ``wiggly'' line) and $\mathrm{MAP}$ point for wavelet prior (green piecewise constant curve with many jumps)}
\end{figure}

A well-studied prior is the Gaussian $N(\mu, \beta(-\Delta)^{-\alpha})$ prior for $\alpha > d-\frac{1}{2}$ and $\beta > 0$ (see, for example, \cite{dashti2013bayesian}). This corresponds to a prior whose samples are randomized weighted Fourier expansions: We say that $v\sim N(\mu,\beta(-\Delta)^{-\alpha})$ in one dimension, if 
\[ v(x) = \mu + \beta \cdot\sum_{k\in\N}\left[ \eps_k \cdot k^{-\alpha} \sin(2\pi k x) + \zeta_k\cdot k^{-\alpha} \cos(2\pi k x)\right],\]
where $\eps_k,\zeta_k\sim N(0,1)$ i.i.d. This is immediate from the Karhunen-Lo\`eve expansion of the Gaussian measure. The two-dimensional analogue is
\begin{align*}
v(x,y) = \mu + \beta \cdot\sum_{k\in\N}&\left[ \eps_k \cdot k^{-\alpha} \sin(2\pi k x)\sin(2\pi l y) + \zeta_k\cdot k^{-\alpha} \cos(2\pi k x)\sin(2\pi l y)\right.\\
&+ \left.  \eta_k \cdot k^{-\alpha} \sin(2\pi k x)\cos(2\pi l y) + \theta_k\cdot k^{-\alpha} \cos(2\pi k x)\cos(2\pi l y)\right]
\end{align*}
with $\eps_k,\zeta_k,\eta_k,\theta_k\sim N(0,1)$ i.i.d. As the trigonometric functions build an orthonormal basis for $L^2([0,1]^d)$, this prior is in principle able to reproduce any permeability function,  but the fact that the covariance operator is diagonal over the Fourier basis leads to a certain type of samples which may be unsuitable in some applications:

Subsurface topology is dominated by ``sheets'' of fairly homogenous permeability with sharp (discontinous) interfaces inbetween (think of a granite layer next to a clay layer). Inference of such a geometry by a Gaussian prior with a Laplacian covariance operator suffers the usual drawbacks of Fourier analysis when applied to ``jumpy'' signals: Sharp interfaces are either not well represented or the Gibbs phenomenon leads to unphysical artefacts around the jump. Also, from a modelling point of view, it does not make a lot of sense to model subsurface topology by periodic base functions. See figure \ref{fig:trigvswavelet} for a demonstration of this phenomenon.

\review{
As the inverse problems consists of ``matching'' (in the Bayesian interpretation of that term) data locally, it is unintuitive to do this with functions which have support on the whole real line, and which do not have any sort of localization. This means that by matching data ``here'', we also change the matching functions values ``there'' (everywhere else). }

An alternative approach is given by priors which are given by randomized weighted series of different base functions than trigonometric polynomials. One possibility is the employment of orthonormal wavelets. This leads to a \textit{localized} orthonormal decomposition of functions. If we use wavelets which are able to model jumps, we can even reproduce sharp interfaces.
\begin{figure}[hbtp]
\centering
\includegraphics[width=0.5\textwidth]{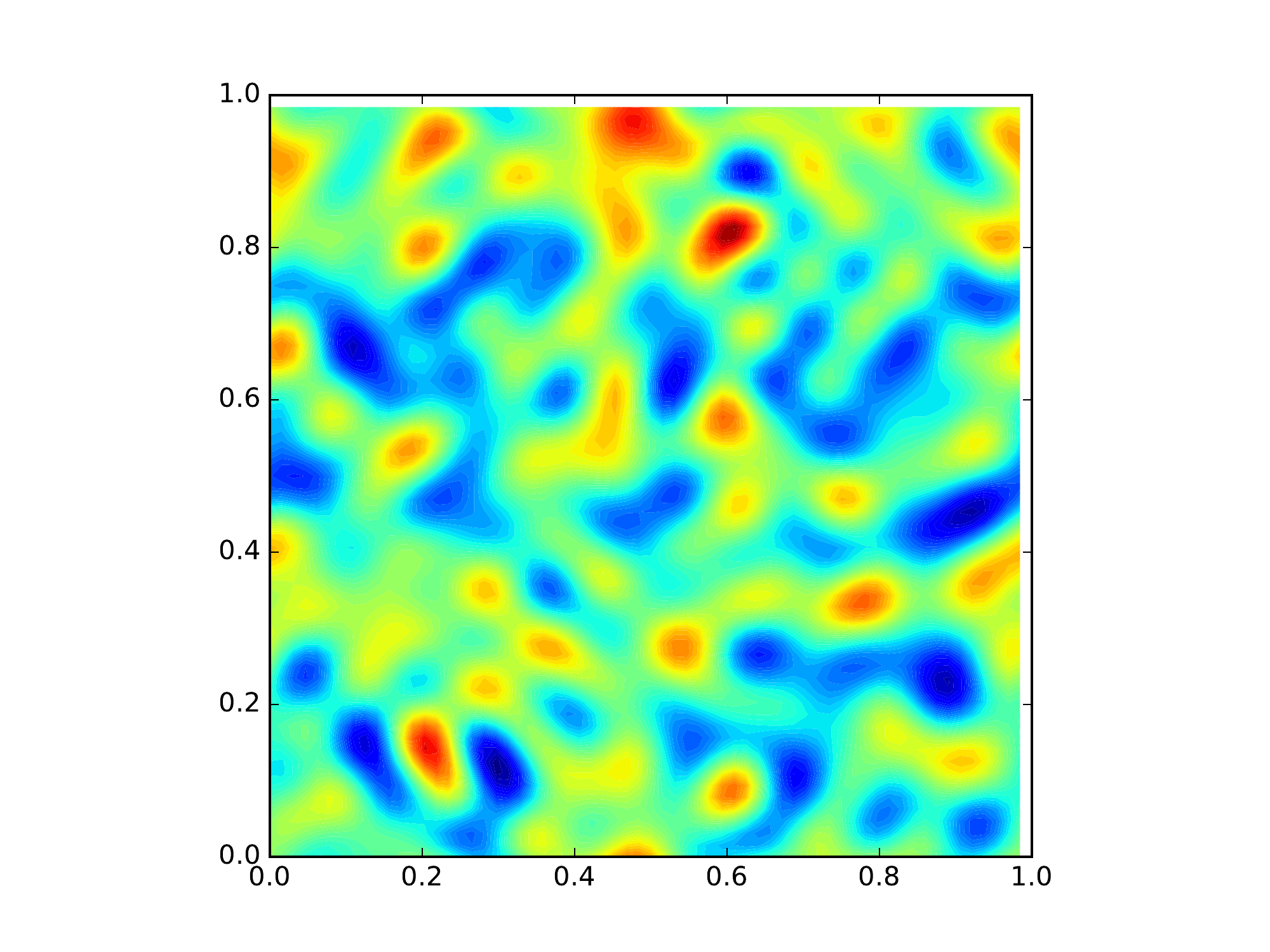}\includegraphics[width=0.5\textwidth]{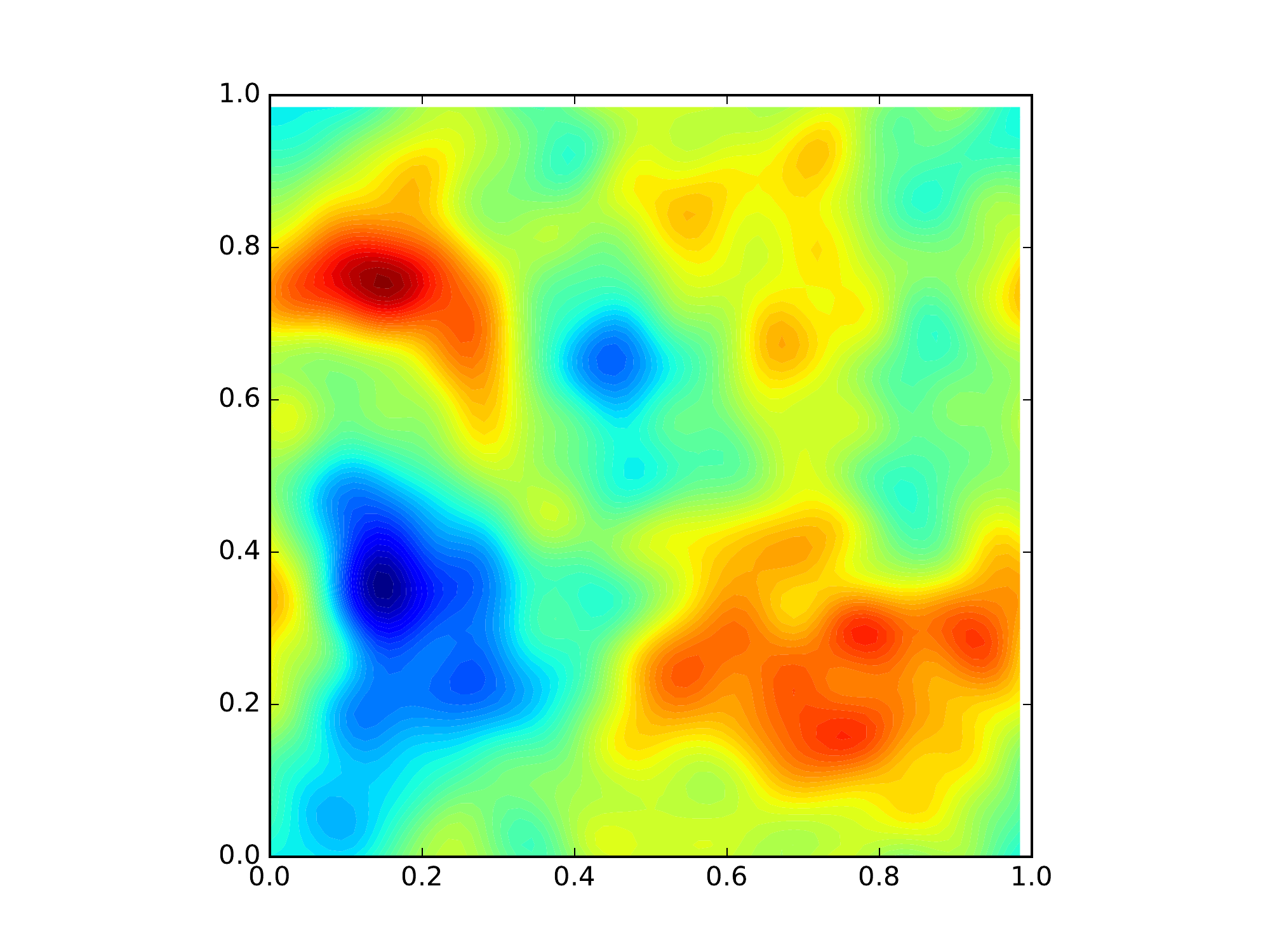}
\caption{Samples from a $N(0,(-\Delta)^{-\alpha})$ prior (left: $\alpha=1.0$, right: $\alpha=2.5$)}
\end{figure}

\review{Use of a wavelet basis in this inverse problem benefits from the same advantages as in the setting of image or audio analysis with wavelets: The localization of the wavelet functions matches nicely the local structure of the observation process

}

\section{Wavelet (``Besov'') priors}\label{sec:waveletpriors}
We can construct a probability measure on Besov spaces\footnote{For a very short introduction to Haar wavelets and Besov spaces consult the appendix. } (i.e. we define random Besov functions) by randomizing coefficients in the wavelet expansion in a suitable way. To this end, we set $\xi_l \sim \exp(-|x|^p/2)$ (i.e. for $p=2$ these are normal random variables, for $p=1$ we have Laplace samples etc.). Furthermore we define an inverse scale parameter $\kappa$ and set $s \geq 0$. Then in dimension $d$ we can define a random function on $[0,1]^d$ by setting
\[ u(x) = \sum_{l=1}^\infty l^{-(\frac{s}{d}+\frac{1}{2}-\frac{1}{p})}\cdot \kappa^{-1/p} \cdot \xi_l \cdot \psi_l(x).\]
This is the notation in \cite{lassas2009discretization}, \cite{dashti2011besov}, which is very suitable for proving various regularity results but not very natural from an implementation point of view. To see this, note that the enumeration of $j,k$ (the resolution and shift parameter) by just one index $l$ and the scaling of the coefficients by $l^{-(\frac{s}{d}+\frac{1}{2}-\frac{1}{p})}$ leads to a spatially inhomogeous span of random functions. We follow the notation in \cite{kolehmainen2012sparsity}, \cite{bui2015scalable}, which stems from the form of equations \eqref{eq:besovdim1} and \eqref{eq:besovdim2} and differentiate by dimension.
\paragraph{1D}
Set $\xi_{j,k}\sim \exp(-|x|^p/2)$ i.i.d. Then 
\[v(x) = \kappa^{-1/p}\cdot \sum_{j=0}^\infty  2^{-j(s+\frac{1}{2}-\frac{1}{p})} \sum_{k=0}^{2^j-1}\xi_{j,k}\cdot \psi_{j,k}(x)\]
is a random function defined on $[0,1]$.
\paragraph{2D}
Set $\xi_{j,k,n}^{(m)}\sim \exp(-|x|^p/2)$ i.i.d. Then 
\[ v(x,y) = \kappa^{-1/p} \cdot \sum_{j=0}^\infty 4^{-j(\frac{s}{2}+\frac{1}{2}-\frac{1}{p})}\sum_{m=0}^2\sum_{k=0}^{2^j-1}\sum_{n=0}^{2^j-1}\xi_{j,k,n}^{(m)} \cdot \psi_{j,k,n}^{(m)}(x,y)\]
is a random function defined on $[0,1]\times[0,1]$.
In either dimension, we call such a random function a sample from a $(\kappa, B_{pp}^s)$ Wavelet prior.

For convenience, we record the following regularity result:

\begin{lemma}[from \cite{lassas2009discretization}, modification with $\kappa$ from \cite{dashti2011besov}]
Let $\kappa>0, p\geq 1, s > 0$ and $v$ a sample from a $(\kappa, B_{pp}^s)$ Wavelet prior as above. The following conditions are equivalent:
\begin{itemize}
\item $\|v\|_{B_{pp}^t} < \infty$ almost surely.
\item $\E \exp(\alpha \|v\|_{B_{pp}^t}^p) < \infty$ for any $\alpha \in (0, \kappa/2)$.
\item $t < s-d/p$.
\end{itemize}
\end{lemma}
Note that in particular, a sample from a $(\kappa, B_{pp}^s)$ Wavelet prior will in general \textit{not} have finite $B_{pp}^s$ norm!

In the case $p=2$, we obtain a Wavelet-based Gaussian prior with mean $0$ and covariance operator $C$ with property
\[C \psi_{j, k, n}^{(m)} = \kappa^{-1}\cdot 4^{-js} \psi_{j, k, n}^{(m)},\]
i.e. we can think of the prior as being a diagonal Gaussian over all Wavelets with exponential decay governed by $s$ and length scale given by $\kappa^{-1}$. In this case the Cameron--Martin space is $E = B_{22}^s = H^s$ (which is the space giving its name to the measure) and with Cameron--Martin norm given by $\|\cdot\|_E = \|\cdot\|_{B_{22}^s}$. In this case, the Cameron--Martin norm generates an inner product by
\[\langle u,v\rangle_E = \sum_{j=0}^\infty 4^{js} \cdot \sum_{k=0}^{2^j-1} w_{j,k}(u) \cdot w_{j,k}(v)\]
in dimension one (where $w_{j,k}(u)$ is the Wavelet decomposition coefficient of $u$ and likewise for $w_{j,k}(v)$). In dimension two,
\[\langle u,v\rangle_E = \sum_{j=0}^\infty 4^{js} \cdot \sum_{m=0}^2\sum_{k=0}^{2^j-1}\sum_{n=0}^{2^j-1} w_{j,k,n}^{(m)}(u) \cdot w_{j,k,n}^{(m)}(v)\]

In the case $p\neq 2$ we define $\|\cdot\|_E = \|\cdot\|_{B_{pp}^s}$ nevertheless. In the ``proper'' Besov case (i.e. non-Gaussian), $E$ is not a Hilbert space (and not the Cameron--Martin space which is only defined for Gaussians). Also, it does not coincide with the space of admissible shifts, which is a Besov space of slightly different regularity. For a proof of this, consult \cite[Lemma 3.5]{agapiou2017sparsity}, see also \cite{helin2015maximum}.

\begin{figure}[hbtp]
\centering
\includegraphics[width=0.5\textwidth]{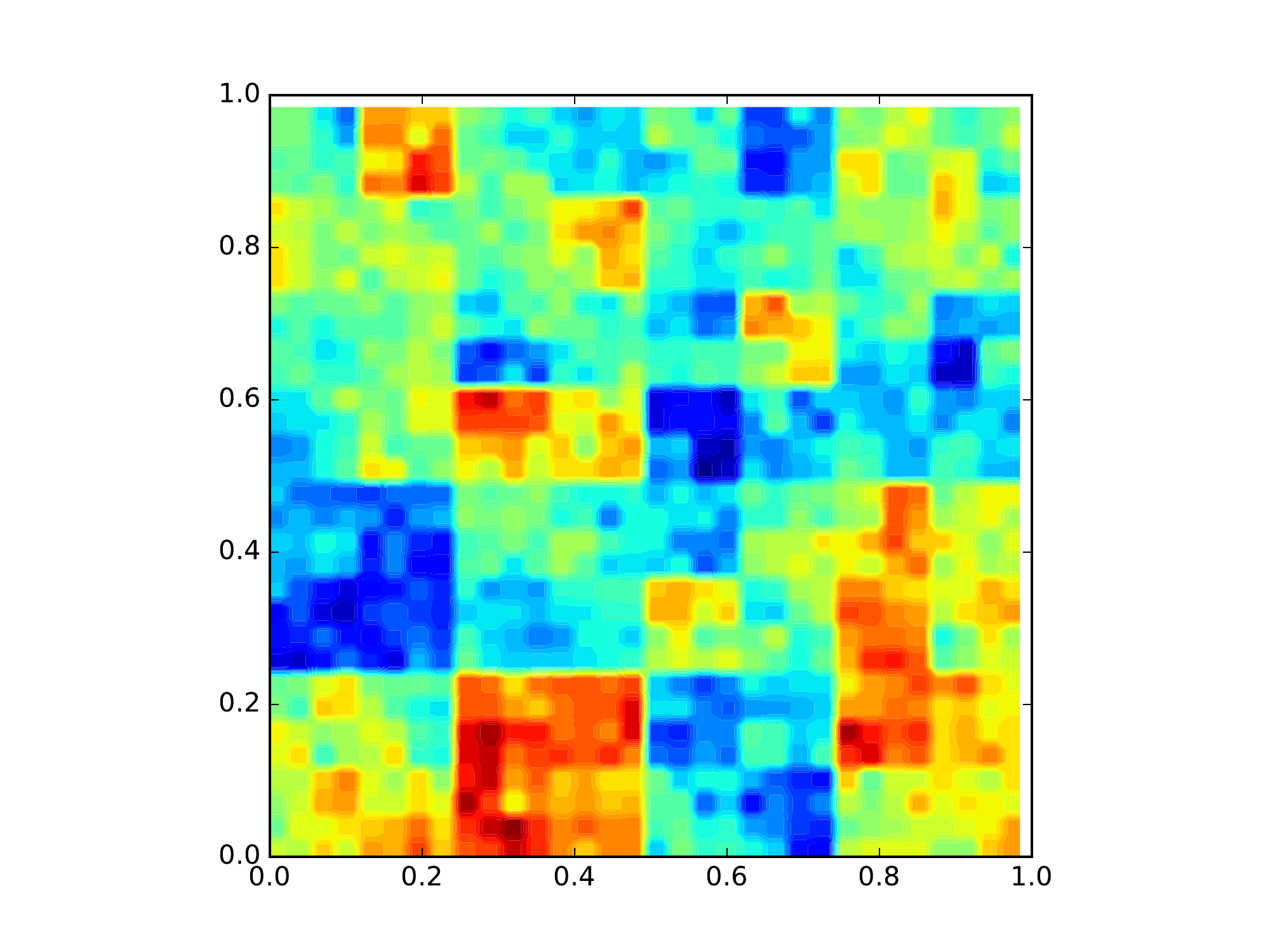}\includegraphics[width=0.5\textwidth]{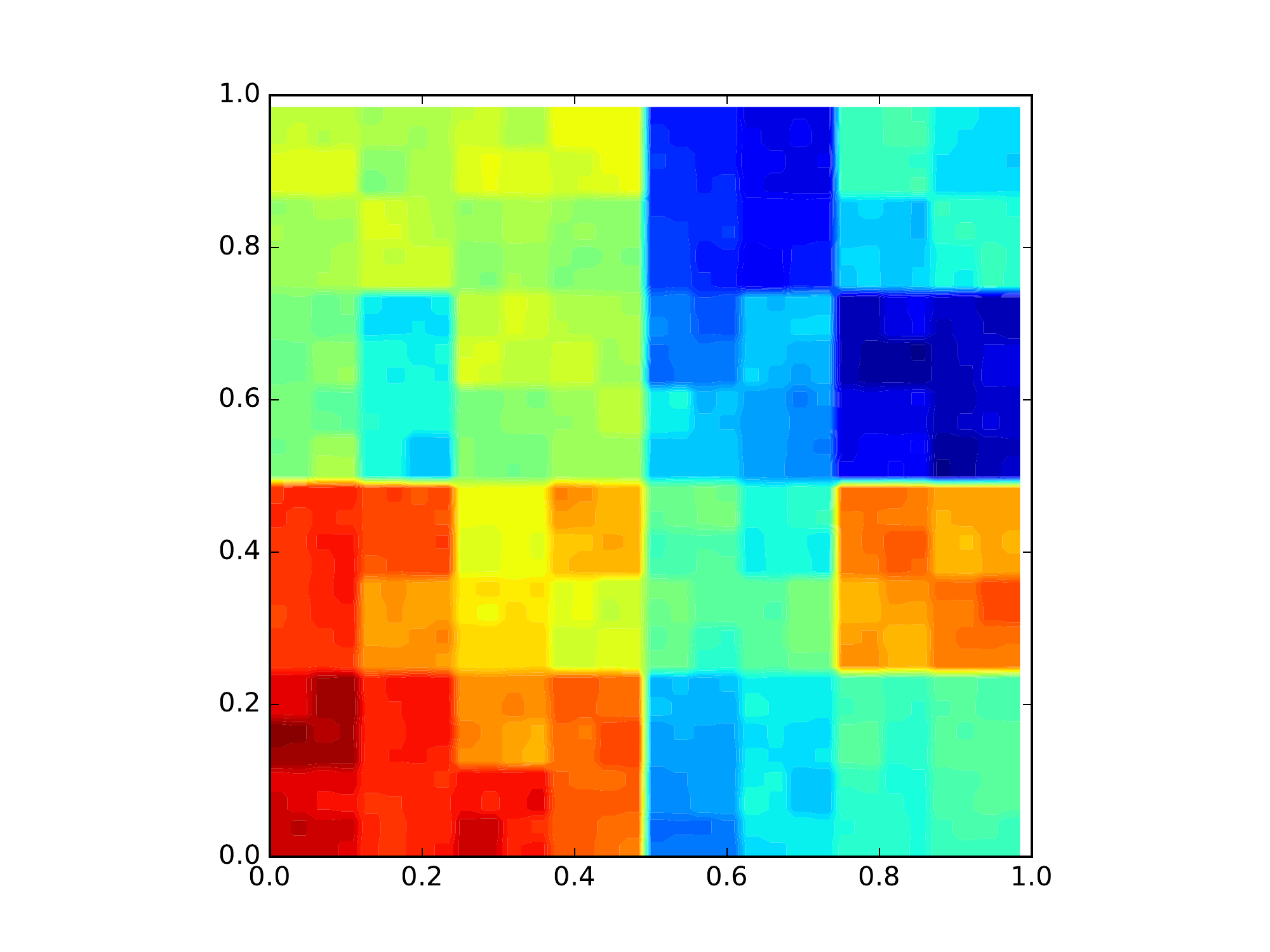}
\caption{Samples from a $(0.1,B_{22}^s)$ prior (left: $s=1.0$, right: $s=2.0$)}
\end{figure}

\subsection{MAP-sparsifying priors} ``Sparse'' wavelet priors (as used in \cite{bui2015scalable}) refer to $B_{1,1}^s$ priors. Here, the MAP point is the optimizer of an energy functional where the prior-dependent part is sparsity-promoting. As pointed out before, this is the only way in which this prior is sparsifying. This is due to the fact that samples from the Laplace distribution (this is the distribution governing the coefficients of the wavelet decomposition) are almost surely non-zero. In this case, the simulations can be carried out in an equivalent way (although the case $p > 1$, $p\approx 1$ is easier to handle as the energy functional is differentiable) and the MAP point is seen to be more sparse.

\subsection{A remark on the name ``Besov prior''}
Note that the notion of randomizing a wavelet expansion in the context of Bayesian inverse problems has been called Besov priors, most notably by \cite{lassas2009discretization}, \cite{dashti2011besov}, \cite{kolehmainen2012sparsity}. In this paper, the name ``Wavelet prior'' is preferred, for several reasons:
\begin{itemize}
\item If the prior is a Gaussian prior (i.e., $p=2$), the $B_{pp}^s$ Besov prior (in the convention of the references above) is named after its Cameron-Martin space and in particular, samples of a $B_{pp}^s$-prior will with probability one \textit{not} lie in $B_{pp}^s$, which can lead to some confusion. 
\item Strictly speaking, in order to construct a $B_{pp}^s$ Besov prior, we need to start with a $\psi$, $\phi$ pair which is at least $s$-regular. A function $f$ is $s$-regular, if $f\in C^s$ and $\partial^\alpha f(x)| \leq C_m(1+|x|)^{-m}$ for all $m\in \N$ and $|\alpha| \leq s$. Our choice of $\phi$ and $\psi$ is not even continuous (but still leads to an interesting prior).
\item From an applied modelling point of view, the fact that ``plausible'' parameters (as defined by the prior) are spanned by a certain class of wavelets (as opposed to trigonometric functions) is more strongly defining feature than its approximation property of some function space.
\end{itemize}
To summarize, in this context, ``Wavelet prior'' means some sort of randomized wavelet expansion which is well-defined in the limit of arbitrarily fine spatial discretization.

\section{Finding the maximum a posteriori point}
Although the construction of the inverse problem was motivated from a Bayesian viewpoint, we will concentrate on obtaining the MAP (maximum a posteriori) reconstruction from the data, as opposed to sampling from the full posterior. This is something that should be done in a subsequent study, although our focus on the MAP estimator can be motivated from the small noise setting: Here we assume that the observations are very close to the truth and thus there is very little variation in samples from the posterior and the MAP estimator is a characteristic quantity. 

The techniques in this section (i.e. using the adjoint approach to calculate quantities of interest) are fairly standard (see for example \cite{giles2006smoking}, \cite{giles2000introduction}). We record them here briefly as in the wavelet setting there is a slight twist which speeds up the computation immensely.
Recall that (write $z=(x,y)$)

\begin{align*}
-\nabla\cdot (e^{u(z)} \cdot \nabla p_u(z)) &= f(z)\text{ in $[0,1]^2$}\\
p(z) &= g(z) \text{ on $\Gamma_D$}\\
\partial_\nu p(z) &= h(z) \text{ on $\Gamma_N$}
\end{align*}

Note that we now record the $u$-dependence of the solution of the PDE explicitly by writing $p_u$.

We first need to find a $G\in H^1(\Omega)$, such that $\operatorname{trace}(G) = g$ and then we define $P_u:=p_u-G$. Then the PDE for $P_u$ is 

\begin{align*}
-\nabla\cdot (e^{u(z)} \cdot \nabla P_u(z)) &= f(z) + \nabla\cdot(e^{u(z)}\cdot\nabla G(z)),\quad\text{ in } [0,1]^2\\
P_u(z) &= 0,\quad \text{ on } \Gamma_D\\
\partial_\nu P_u(z) &= h(z)-\partial_\nu G(z),\quad \text{ on } \Gamma_N.
\end{align*}

A good space of test functions is $H_{0,D}^1 := \{\phi\in H^1(\Omega):\, \operatorname{trace}(\phi) = 0 \text{ on } \Gamma_D\}$. Now this has a variational form: We call $P_u\in H_{0,D}^1$ the solution of  the PDE, if $P_u=P$ is the unique solution of

\[ \int_\Omega e^u\cdot  \nabla P\cdot \nabla \phi = \int_\Omega f\phi - \int_\Omega e^u\cdot \nabla\phi\cdot \nabla G + \int_{\Gamma_N} h \operatorname{trace}(\phi)
\quad \text{for all } \phi\in  H_{0,D}^1.\]

Then $p_u = P_u + G$ is the solution of the original PDE. We can write that $p_u \in  H_{g,D}^1 :=\{\phi\in H^1(\Omega):\, \operatorname{trace}(\phi) = g \text{ on } \Gamma_D\} $

\subsection{\Frechet-derivative of a Quantity of Interest}

We assume to have a quantity of interest\footnote{This section closely follows \cite{sullivan2015introduction}. Thanks to Donald Estep for explaining to me how to obtain the adjoint PDE with inhomogenous boundary conditions, see also \cite{estep2004short}.} depending on the solution to our original PDE, $q:H_{g,D}^1\to\R$, and we are interested in the sensitivity of $q(p_u)$ on the parameter scalar field $u$. 

For this to work we define the implicit dependence of $u$ and $p_u$ by the following expression: $p_u$ is the unique (by PDE theory) $p\in H_{g,D}^1$ such that 
\[ 0 = J(u, p) = -\nabla\cdot(e^u\cdot \nabla p) - f.\]
The implicit function theorem yields the existence of a function $\mathfrak p: L^\infty \to H_{g,D}^1$ with $J(u, \mathfrak p(u)) = 0$ and $D_u \mathfrak{p}(\bar u) = [D_p J(\bar u, \mathfrak p(\bar u))]^{-1}\circ D_u J(\bar u, \mathfrak p(\bar u))$.
Then our quantity of interest can be written as a function of $u$.\footnote{We will use $p_u$ and $\mathfrak{p}(u)$ interchangeably depending on whether we view it as a $H_{g,D}^1$ object or whether we focus on the explicit $u$-dependence coming from the implicit function theorem}
\[ Q(u) := q(p_u) = q(\mathfrak{p}(u)).\]
We can calculate the mentioned sensitivity by evaluating the \Frechet~derivative
\begin{align}
D_u Q(\bar u) &= D_p q(p_{\bar u}) \circ D_u \mathfrak p(\bar u) =  D_p q(p_{\bar u}) \circ [D_p J(\bar u, \mathfrak p(\bar u))]^{-1}\circ D_u J(\bar u, \mathfrak p(\bar u))\\
D_u Q(\bar u) [\varphi]&= D_p q(p_{\bar u})\left\{ [D_p J(\bar u,  p_{\bar u})]^{-1}\left( D_u J(\bar u, p_{\bar u})[\varphi]\right)\right\}\label{eq:sensitivityQ}
\end{align}
Now we use the special structure of $J$:
\begin{align*}
D_u J(\bar u, p_{\bar u})[\varphi] &= -\nabla \cdot \left(\varphi\cdot e^{\bar u}\cdot \nabla p_{\bar u}\right)\\
D_p J(\bar u, p_{\bar u})[v] &= -\nabla \cdot \left( e^{\bar u} \cdot \nabla v\right)
\intertext{and hence}
w &= [D_p J(\bar u,  p_{\bar u})]^{-1}\left( D_u J(\bar u, p_{\bar u})[\varphi]\right)
\end{align*}
is the unique element $w$ given by the PDE
\begin{equation}\label{eq:primal}
\begin{split}
-\nabla \cdot \left( e^{\bar u} \cdot \nabla w\right) &= -\nabla \cdot \left(\varphi\cdot e^{\bar u}\cdot \nabla p_{\bar u}\right)\\
w &= 0\text{ on $\Gamma_D$}\\
\partial_\nu w &= 0\text{ on $\Gamma_N$}.
\end{split}
\end{equation}
The boundary conditions on $w$ arise from the fact that $D_u J(\bar\theta, u_{\bar\theta}) : H_{0,D}^1 \to \R$ because $H_{0,D}^1$ is the ``difference space'' of $H_{g,D}^1$.
Combining this with \eqref{eq:sensitivityQ}, we can express $D_u Q(\bar u)[\varphi]$ as follows:
\[ D_u Q(\bar u)[\varphi]= D_p q(p_{\bar u})[w], \]
which is a linear functional acting on $w$. Now if we are interested in the gradient of $Q(\bar\theta)$, there is a straightforward way to do this by calculating $ D_\theta Q(\bar \theta)[\varphi]$ for ``all'' directions $\varphi$ and concluding this in a column vector which then is the gradient. This entails solving the PDE on $w$ as many times as we have directions $\varphi$, which can be prohibitively costly. By the adjoint method, we have a more clever way.
\subsubsection{Adjoint functional evaluation}
We consider a second-order elliptical PDE of the form
\begin{align*}
-\nabla \cdot (e^u\cdot \nabla v) &= f \text{ in } \Omega\\
v &= g \text{ on } \Gamma_D\\
\partial_\nu v &= h \text{ on } \Gamma_N
\end{align*}
for $f\in H^{-1}$, $\Gamma_D\cup\Gamma_N = \partial \Omega$ and we are interested in the value of 
\[\int_\Omega \ell\cdot v,\]
where $\ell\in H^{-1}$. Now this is equivalent to evaluating \[\int_\Omega f\cdot\tilde v - \int_{\Gamma_D}g \cdot e^\theta\cdot  \partial_\nu \tilde v + \int_{\Gamma_N} h\cdot e^\theta\cdot \tilde v,\]
where 
\begin{align*}
-\nabla \cdot (e^u\cdot \nabla \tilde v) &= l \text{ in } \Omega\\
\tilde v &= 0 \text{ on } \Gamma_D\\
\partial_\nu  \tilde v &= 0 \text{ on } \Gamma_N
\end{align*}
as can be shown easily by partially integrating twice.
\subsection{\Frechet-derivative of a Quantity of Interest (continued)}
As \[  D_u Q(\bar u)[\varphi]= D_p q(p_{\bar u})[w]\]
is a linear functional acting on $w$, we are exactly in the setting of the last section. 

We know that we can evaluate this quantity by working with the adjoint PDE for $w$ which is
\begin{equation}\label{eq:dual}
\begin{split}
-\nabla \cdot \left( e^{\bar u} \cdot \nabla \tilde w\right) &= D_p q(p_{\bar u})\text{ in $[0,1]^2$}\\
\tilde w &= 0\text{ on $\Gamma_D$}\\
\partial_\nu \tilde w &= 0\text{ on $\Gamma_N$}.
\end{split}
\end{equation}
This is just one PDE which we need to solve and now we just need to evaluate the inner product of $\tilde w$ with the primal PDE's \eqref{eq:primal} right hand side
\[\int_\Omega -\nabla \cdot (\varphi e^{\bar u} \cdot \nabla p_{\bar u})\cdot \tilde w = \int_\Omega \varphi\cdot e^{\bar u}\cdot \nabla p_{\bar u}\cdot \nabla \tilde w\]
for each direction $\varphi$ and this will be $ D_u Q(\bar u)[\varphi]$. In the case that our parameter space is spanned by a Wavelet basis and we are interested in the sensitivity in direction $\phi$, where $\phi$ is an element of our wavelet basis, the inner product is actually the computation of the wavelet decomposition of $e^{\bar u}\cdot \nabla p_{\bar u} \cdot \nabla \tilde w$: 

\begin{align*}
D_u Q(\bar u)[\psi_l] &= c_l, \text{ where}\\
e^{\bar u}\cdot \nabla p_{\bar u} \cdot \nabla \tilde w&= \sum_{l=1}^\infty c_l\cdot \psi_l 
\end{align*}
is a Wavelet decomposition and $\tilde w$ is the solution of \eqref{eq:dual}.

For a wavelet resolution of $2^J\times 2^J$ we have $4^J$ ``directions'', i.e. $4^J$ entries in the gradient. This means we have diminished our computational effort for computing the gradient from  $4^J$ PDE solves to $2$ PDE solves (one solve for $p_{\bar u}$ and one solve for $\tilde w$ which depends on $p_{\bar u}$) and one wavelet decomposition calculation, which in general is a huge improvement. 
\subsection{Sensitivity of a data misfit functional}\label{sec:sensitivityDMF}
Consider the case of \[q(p) = \Phi(p) = \frac{1}{2\gamma^2}\|y-\Pi p\|^2\] in the setting of section \ref{sec:model}. Its \Frechet\, derivative is $D_p q(p_{\bar u})(w) = -\frac{1}{\gamma^2}\cdot \langle y-\Pi p_{\bar u}, \Pi w\rangle$ and we can write $D_p q(p_{\bar u})\in H^{-1}$ by 
\[D_p q(p_{\bar u}) = -\frac{1}{\gamma^2}\cdot \sum_{i=1}^N (y_i - p_{\bar u}(x_i))\cdot \delta_{x_i}\]
and hence we can compute the sensitivity w.r.t. the parameter $u$ by the following two-step procedure:
\begin{algo}[Computing the QoI's gradient]~\\

\begin{enumerate}
\item Compute the solution $\tilde w$ of 
\begin{align*}
-\nabla \cdot \left( e^{\bar u} \cdot \nabla \tilde w\right) &=-\frac{1}{\gamma^2}\cdot \sum_{i=1}^N (y_i - p_{\bar u}(x_i))\cdot \delta_{x_i}\\
\tilde w &= 0\text{ on $\Gamma_D$}\\
\partial_\nu \tilde w &= 0\text{ on $\Gamma_N$}.
\end{align*}
\item Obtain the wavelet expansion of $e^{\bar u}\cdot \nabla p_{\bar u}\cdot \nabla \tilde w$ in the form
\[ e^{\bar\theta}\cdot \nabla u_{\bar\theta}\cdot \nabla \tilde w = \sum_{l=1}^\infty c_l \cdot \psi_l.\]
There are two ways to do that.\begin{enumerate}
\item \review{(Method 1): Compute $c_l = \int_\Omega \psi_l\cdot e^{\bar u}\cdot \nabla p_{\bar u}\cdot \nabla \tilde w$ }
\item \review{(Method 2): Calculate the fast wavelet transformation on a grid.}
\end{enumerate}
Then $D_u Q(\bar u)(\psi_l) = c_l$, or 
\[ \nabla_u Q(\bar u) = (c_l)_{l\in\N}.\]
\end{enumerate}
\end{algo}
Note that the last step is the only step where wavelets come explicitly into play. The equivalent version for a trigonometric prior would be to look for the Fourier decomposition instead, method 2 being the calculation of the Fast Fourier transform. There is a major difference, though:
\begin{rem}\label{rem:main}
The last step is where the choice of wavelets really pays off: The solution $\tilde w$ of the dual PDE encodes the gradient information and entries of the gradient vector can be calculated by computing the projection of $e^{\bar u}\cdot \nabla p_{\bar u}\cdot \nabla \tilde w$ onto the elements of the basis used (in this case: wavelet functions). This is easily approximated for the choice of a wavelet basis by constructing the wavelet decomposition of said function. If we take trigonometric functions, this approximation is a lot worse due to \review{two reasons: First, the non-periodicity of $e^{\bar u}\cdot \nabla p_{\bar u}\cdot \nabla \tilde w$: This means that the Fourier-based approximation will have bad approximative properties at the boundaries of the domain. Second, wavelets are able to approximate the quantity pointwise and arbitrarily well (in the $\|\cdot\|_\infty$-norm), which is a crucial property in this setting. A Fourier-based approximation will fail to do that given a sufficiently irregular adjoint solution (note that the quantity $e^{\bar u}\cdot \nabla p_{\bar u}\cdot \nabla \tilde w$ is theoretically only guaranteed to be in $L^2$).}
\end{rem}
\subsection{Computing the MAP point for the $B_2^s$ prior}
Now that everything is in place, computing the MAP point consists of solving the optimization problem 
\[u_{\mathrm{MAP}} = \operatorname{arginf} I(u) =\operatorname{arginf} \Phi(u) + \frac{1}{2}\|u\|_E^2\]
with $\|\cdot\|_E$ defined as in section \ref{sec:waveletpriors}. We will concentrate on the Gaussian case $p=2$ where the role of $E$ as the Cameron--Martin space of the prior is well-known.

This can be done with a standard BFGS method, for which we need the gradient of $I$. After the last section, this is now immediately given by (note that $h$ is a general direction in which the gradient is evaluated)
\[ DI(u)[h] = D\Phi(\mathfrak{p}(u))[h] + \langle u, h\rangle_E\]
where $D\Phi(\mathfrak{p}(u))[h] = D_u Q(u)[h]$ is calculated as described in section \ref{sec:sensitivityDMF}.
\review{
\subsection{Computing the MAP point for the $B_1^s$ prior}
Similarly, the MAP point for the Bayesian inverse problem with the "sparse" $B_1^s$ prior is given by
\[u_{\mathrm{MAP}} = \operatorname{arginf} I(u) =\operatorname{arginf} \Phi(u) + \|u\|_{B_1^s}.\]
We formulate this in terms of $u$'s wavelet decomposition. Let $W$ be the inverse wavelet transform operator, i.e. $W$ maps a wavelet decomposition into function space, then 
\[w_{\mathrm{MAP}} = \operatorname{arginf} I(Ww) + |w_0| + \sum_{j=0}^\infty 2^{j(s-d/2)}\sum_{m=0}^2\sum_{k=0}^{2^j-1}\sum_{n=0}^{2^j-1}|w_{j,k,n}^{(m)}|\]
and 
\[u_{\mathrm{MAP}} = Ww_{\mathrm{MAP}}.\]
This minimization problem is readily solved with the FISTA algorithm \cite{beck2009fast}, again benefitting from the particularly fast gradient calculation for a Wavelet-based prior.
}
\section{Simulations}
\review{
All simulations refer to the elliptical inverse problem from section \ref{sec:model}. We consider fixed (and known) boundary conditions and source terms and in particular various ground truth permeabilities. Note that in the following, $J_{\text{max}}$ is the highest resolution of the wavelet decompositions (e.g. for the Besov prior) used, i.e. any function in wavelet decomposition is constructed as 
\[ v(x,y) = \kappa^{-1/p} \cdot \sum_{j=0}^{J_{\text{max}}} 4^{-j(\frac{t}{2}+\frac{1}{2}-\frac{1}{p})}\sum_{m=0}^2\sum_{k=0}^{2^j-1}\sum_{n=0}^{2^j-1}\xi_{j,k,n}^{(m)} \cdot \psi_{j,k,n}^{(m)}(x,y).\]

Note that we specifically use the Haar wavelet basis here.

To avoid an inverse crime (\cite{kaipio2006statistical}), all grids used for numerical solution of the PDE were at least four to eight times finer than the maximal resolution $J_{\text{max}}$ used by the wavelet decomposition. Also, observation data was not collected at grid points but on points (randomly chosen or set on a fixed grid which is not the grid generated by the PDE solver) using interpolation of the evaluation function.
All PDEs (primal and adjoint) are solved with \texttt{FEniCS} on the domain using a regular triangular mesh (with mesh size at least four times finer than $2^{J_{\text{max}}}$). 

To demonstrate how much more efficiently gradient information can be obtained in the setting of a wavelet-based prior, we compare four different calculations of the MAP estimator:

\begin{itemize}
\item Method 1 for a trigonometric prior
\item Method 2 for a trigonometric prior
\item Method 1 for a wavelet-based prior
\item Method 2 for a wavelet-based prior
\end{itemize}

To simplify the exposition and because we are for now only interested in the relative computational advantages of using a wavelet-based prior in contrast to a trigonometric prior, we compare all methods in the same way:

We fix an artificial ground truth permeability which is not a sample of either prior but can be thought of as being a semi-realistic toy model. Then, for a given set of noisy observations on a grid of observation points, we try to obtain the MAP estimators, respectively. 

The ground truth logpermeability is given by the following function:
\begin{align*} u_{\text{truth}}(x, y) &= 0.5/\log_{10}(e)\cdot \sin(4\cdot x) +\\
	&+ \begin{cases}-2/\log_{10}(e) & \text{ if }  x \leq 0.85 \text{ and } (x-1)^2+y^2 \in (0.55^2, 0.65^2)\\
-4/\log_{10}(e) & \text{ if } x \in (0.4375, 0.5) \text{ and } y > 0.625 \end{cases}\end{align*}

The boundary conditions are mixed: Dirichlet-0 condition on $x=0$ and Neumann-0 condition on the remaing three edges.

There are three circular sources, yielding a source function
\[ f(x,y) = \chi_{B_{0.1}(0.6, 0.85)}\cdot (-2000) + \chi_{B_{0.1}(0.2, 0.75)}\cdot (2000) + \chi_{B_{0.1}(0.8, 0.2)}\cdot (-2000).\]

\begin{figure}[hbtp]
\centering
\includegraphics[width=0.5\textwidth]{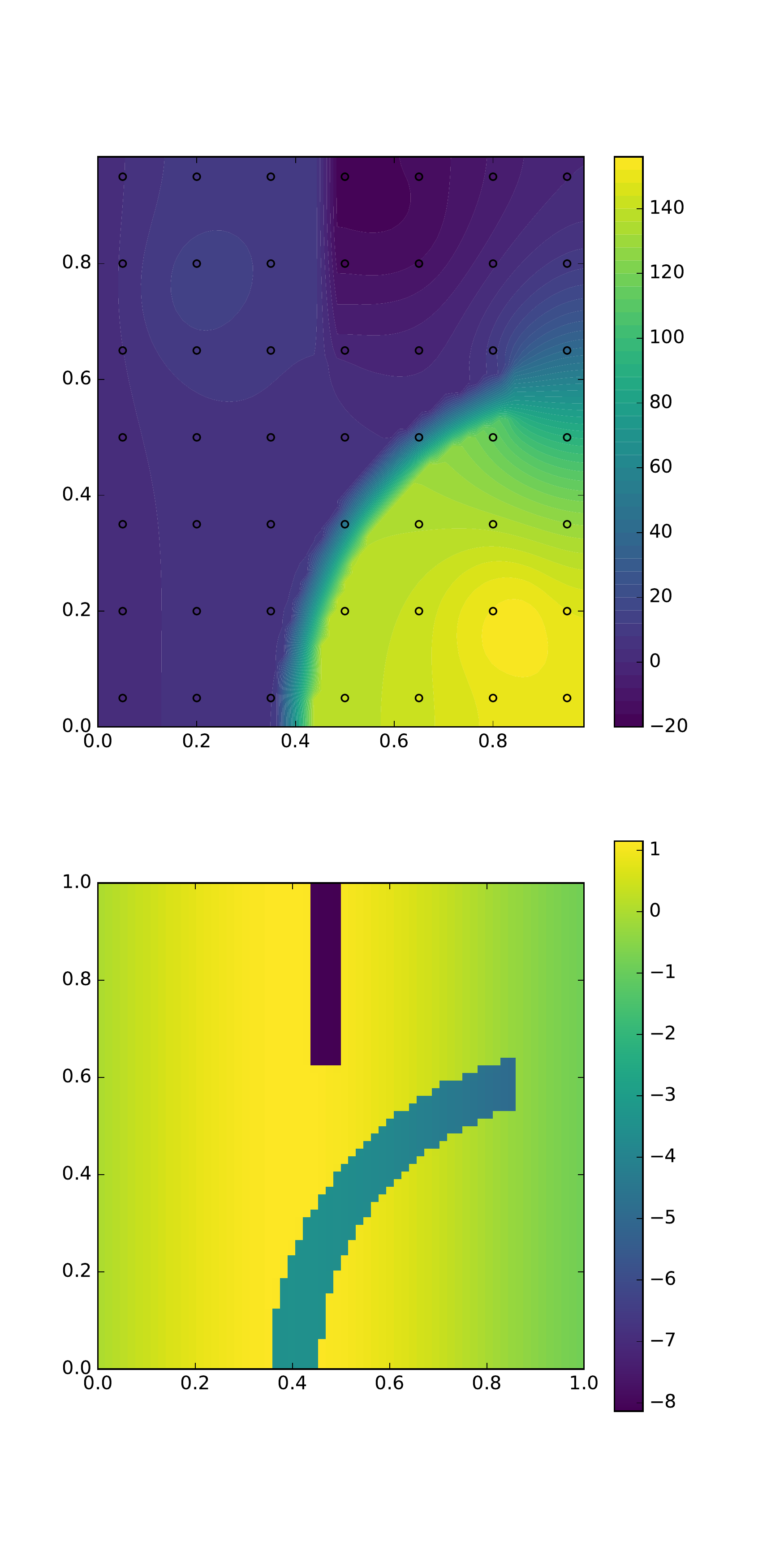}
\caption{The ground truth permeability (bottom) and resulting pressure field with pointwise measurement locations (top).}
\end{figure}

The observational positions are a $7\times 7$ regular grid with vertices \[(0.05, 0.05), (0.05, 0.95), (0.95, 0.05), (0.95, 0.95),\] which is not aligned with or coinciding with any points on the PDE solver's mesh grid, thereby preventing an inverse crime (\cite{kaipio2006statistical}). The observational noise was defined to be Gaussian with mean $0$ and standard deviation $1$.

The PDE is solved on a rectangular grid of size $2^7\times 2^7$ (i.e. 16384 cells).

The priors were specified as follows: The wavelet-based prior is a Gaussian Wavelet $B_2^{1.5}$ prior and the Fourier-based prior measure is the Gaussian trigonometric prior $N(0, 2\cdot(-\Delta)^{-1/2})$ (this is not a proper prior in the discretization limit but here we are not concerning ourselves with discretization invariance anyway). The two priors could in principle be calibrated such that they assign similar ``Cameron-Martin weight'' to the ground truth parameter, but again, we are only concerned with the speed of inference, not the comparability of the priors.
The discretization cutoff of the priors is chosen such that the maximal number of parameters for the trigonometric prior is $33\times 33 = 1089$ and the maximal number of paramaters for the Wavelet-bases prior is $1024$, so that the degrees of freedom are quite comparable. 

Note that the data size is smaller than the number of parameters by a factor of more than 20, with the numerical discretization again being roughly a factor of 16 times that.

Recall that we want to compare four different gradient calculations (and in particular, their varying speed): With a Wavelet prior versus with a Fourier prior, using the two methods respectively(method 1 being the slightly more exact, but slower computational method of obtaining the gradient). It is little revealing to talk about the millisecond advantages of gradient evaluation, so we will instead look at how this affects a real MAP computation. 
For this, we use a simple gradient descent method with a standard backstepping rule:
	\[u_{k+1} = u_k - \frac{\alpha}{2^N}\nabla I(u_k),\]
where $N$ is the smallest integer such that $I(u_{k+1}) < I(u_k) - \alpha/2\|\nabla I(u_k)\|$. The only way all four different calculations differ is the method of how the gradient $\nabla I$ (or rather the difficult part $\nabla \Phi$) is obtained. We do not only track the number of iterations but also consider the time spent in each iteration to compare the efficacy and efficiency of each method. We let the simulation run for roughly $6000s$ each and compare.

These are the findings of the simulation: 
\begin{itemize}
\item The ``quicker'' method 2 does not work at all in the Fourier prior: The ``gradient'' returned by this method does not correspond to any sensible descent direction (the reasons for that are explained in remark \ref{rem:main}) after a few iterations and the optimization procedure terminates unsuccessfully. For this reason, we do not incorporate this option in the figures after that.
\item The ``slower'' method 1 in the Fourier setting manages to reduce the functional $I$ by quite some extent and finishes (by time-out, not by achieving its optimization goal which is not to be expected with such a crude method like gradient descent anyway) with a nevertheless acceptable proposal for a MAP estimator, satisfactorily matching the data.
\item The ``slower'' method 1 in the Wavelet setting is slightly faster than its trigonometric counterpart (but note that the absolute magnitudes of $I$ are not comparable due to the two priors being not "on equal terms" with respect to the artificial ground truth) and produces a similarly acceptable proposal.
\item The ``quicker'' method 2 for the Wavelet prior is so quick that during the course of $6000s$ it manages to converge (as can be seen arguably by observing the levelling out of the functional $I$) to the true MAP estimator. The error made with this slightly inexact calculation of the gradient is unimportant as can be seen in figure \ref{fig:comparison2}, where the cost functional $I$ is plotted not against time, but against iteration count: Here, method 1 and 2 in the Wavelet setting are almost on top of each other, with method 2 being much quicker (as can be seen in figure \ref{fig:comparison} where we incorporate computation time).
\item Not particularly relevant to comparison of the two methods but still interesting is the observation that the not-yet-converged result of method 1 in the wavelet setting shows much finer structure than the eventual MAP estimator yielded by method 2. This phenomenon is mirrored by figure \ref{fig:norms}: The optimization procedure drives the Cameron-Martin-norm of the iteration (penalizing fine structure) first up and is only reduced again after a reasonable match with the data was obtained.
\end{itemize}

\begin{figure}[hbtp]
\centering
\label{fig:comparison}
\begin{minipage}[t]{0.45\linewidth}
\includegraphics[width=\textwidth]{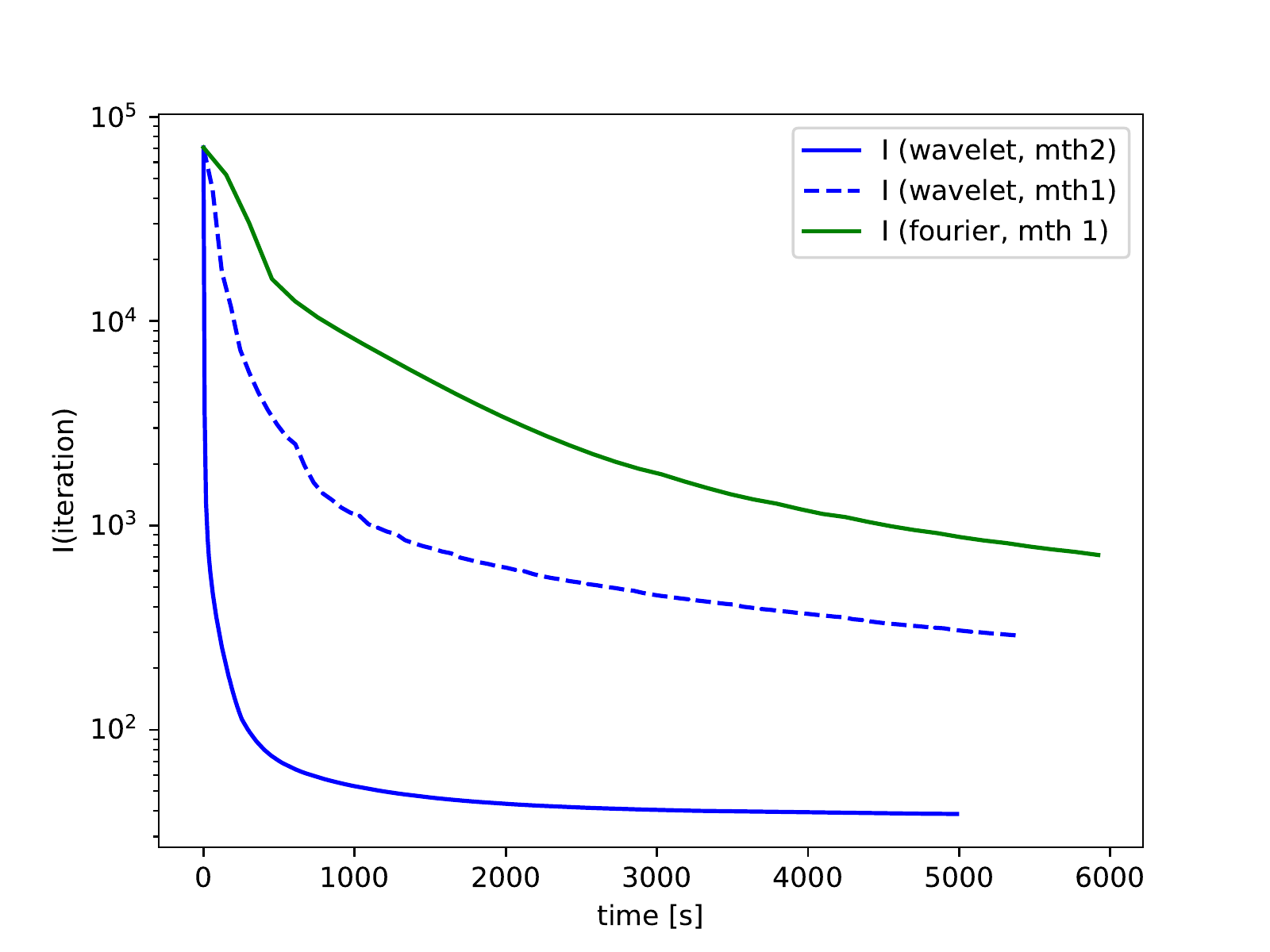}
\caption{Plot of error functional $I$ over computation time in seconds for the three relevant options. Note that only method 2 in the wavelet setting achieves optimality.}
\end{minipage}
\hspace*{0.05\linewidth}
\begin{minipage}[t]{0.45\linewidth}
\includegraphics[width=\textwidth]{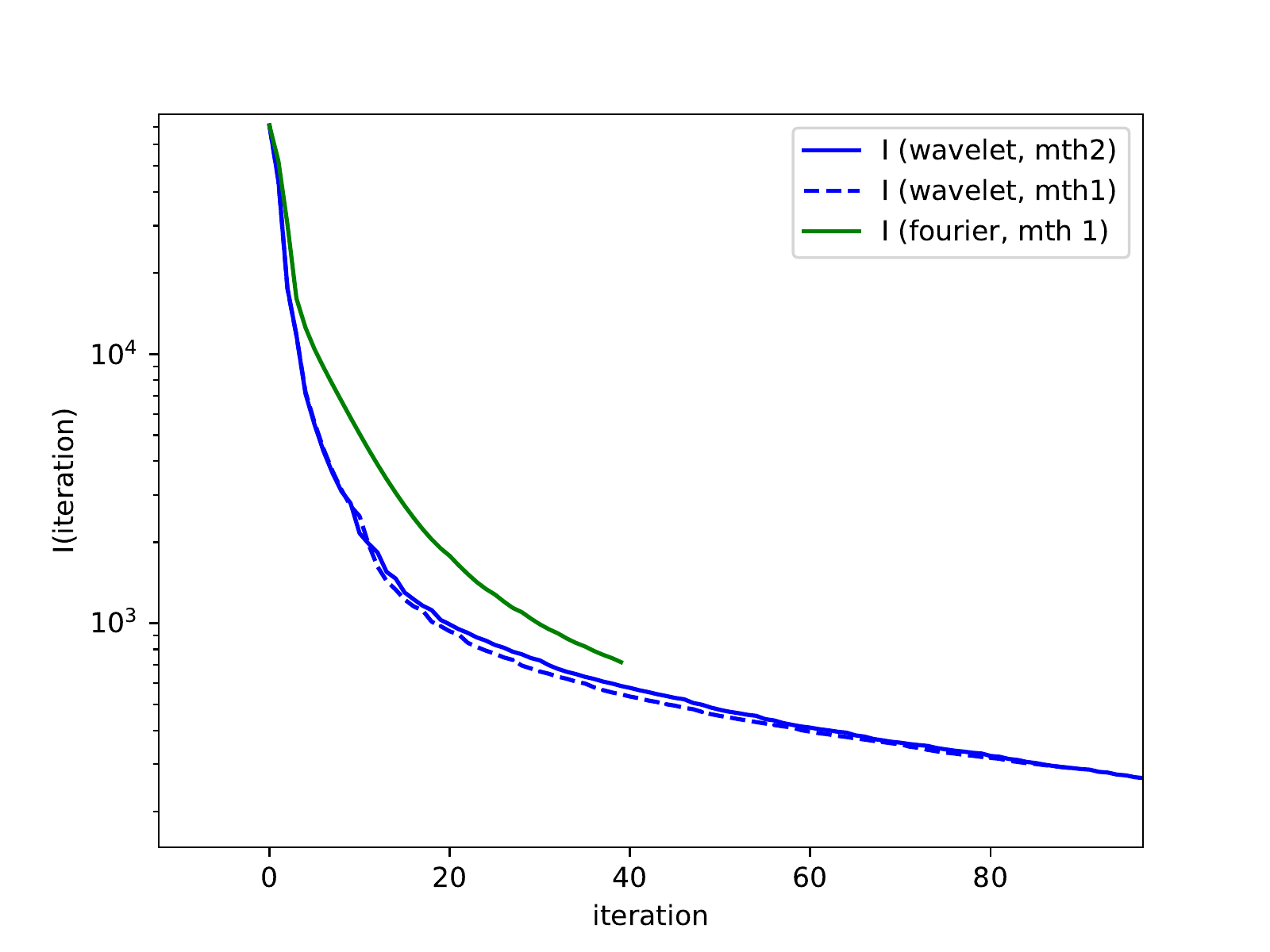}
\caption{Plot of error functional $I$ over iteration number for the three relevant options. Note that both methods in the wavelet setting are essentially equivalent when considered ``per iteration''.}
\label{fig:comparison2}
\end{minipage}
\end{figure}

\begin{figure}[hbtp]
\centering
\begin{minipage}[b]{0.45\linewidth}
\includegraphics[width=\textwidth]{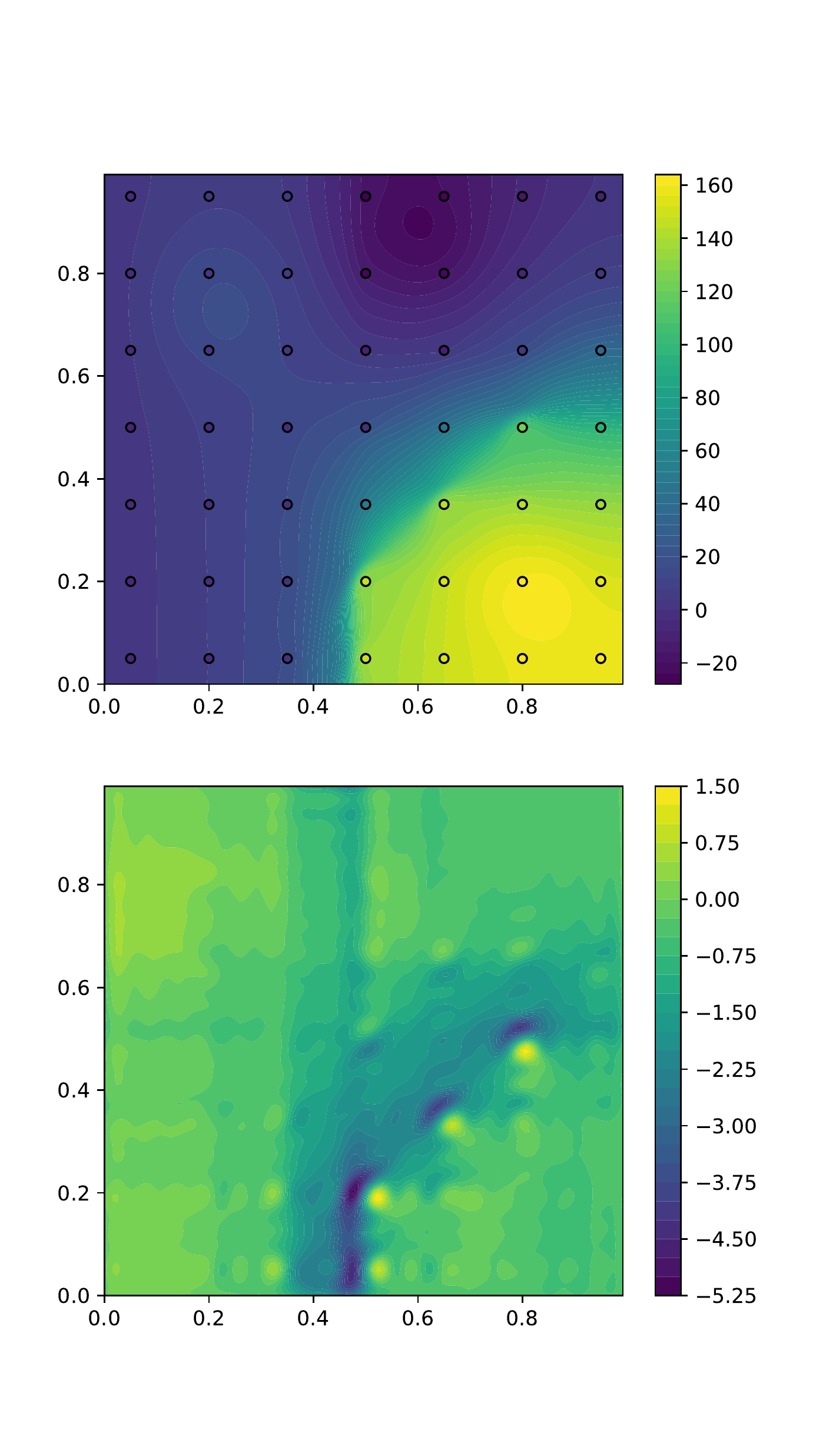}
\caption{MAP ``proposal'' (i.e. not converged result) for method 1 in the Fourier setting}
\end{minipage}
\hspace*{0.05\linewidth}
\begin{minipage}[b]{0.45\linewidth}
\includegraphics[width=\textwidth]{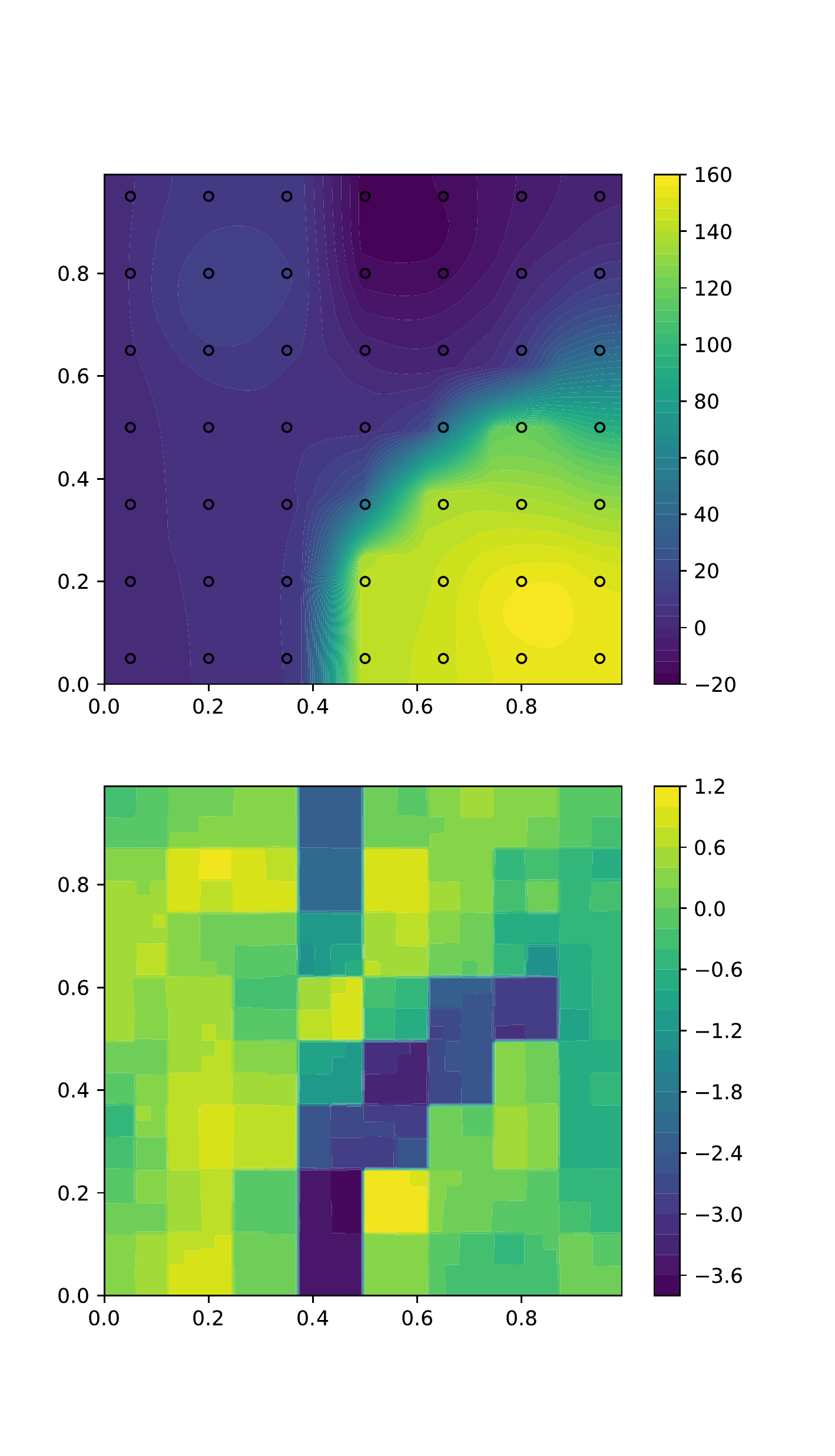}
\caption{MAP estimator (converged) for method 2 in the Wavelet setting}
\label{fig:wavelet_MAP}
\end{minipage}
\begin{minipage}[b]{0.45\linewidth}
\includegraphics[width=\textwidth]{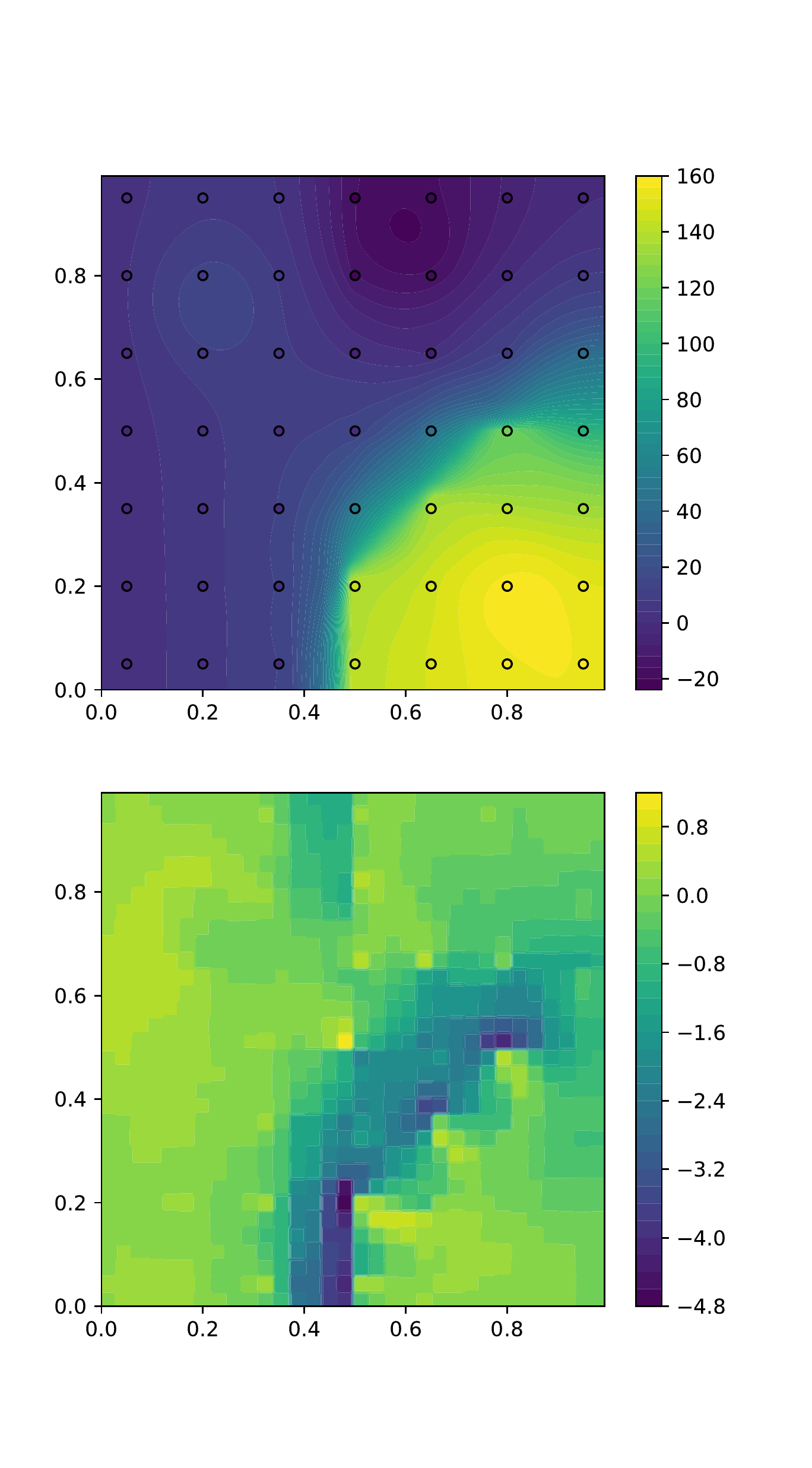}
\caption{MAP ``proposal'' (i.e. not converged result) for method 1 in the wavelet setting}
\label{fig:wavelet_proposal}
\end{minipage}
\hspace*{0.05\linewidth}
\begin{minipage}[b]{0.45\linewidth}
\includegraphics[width=\textwidth]{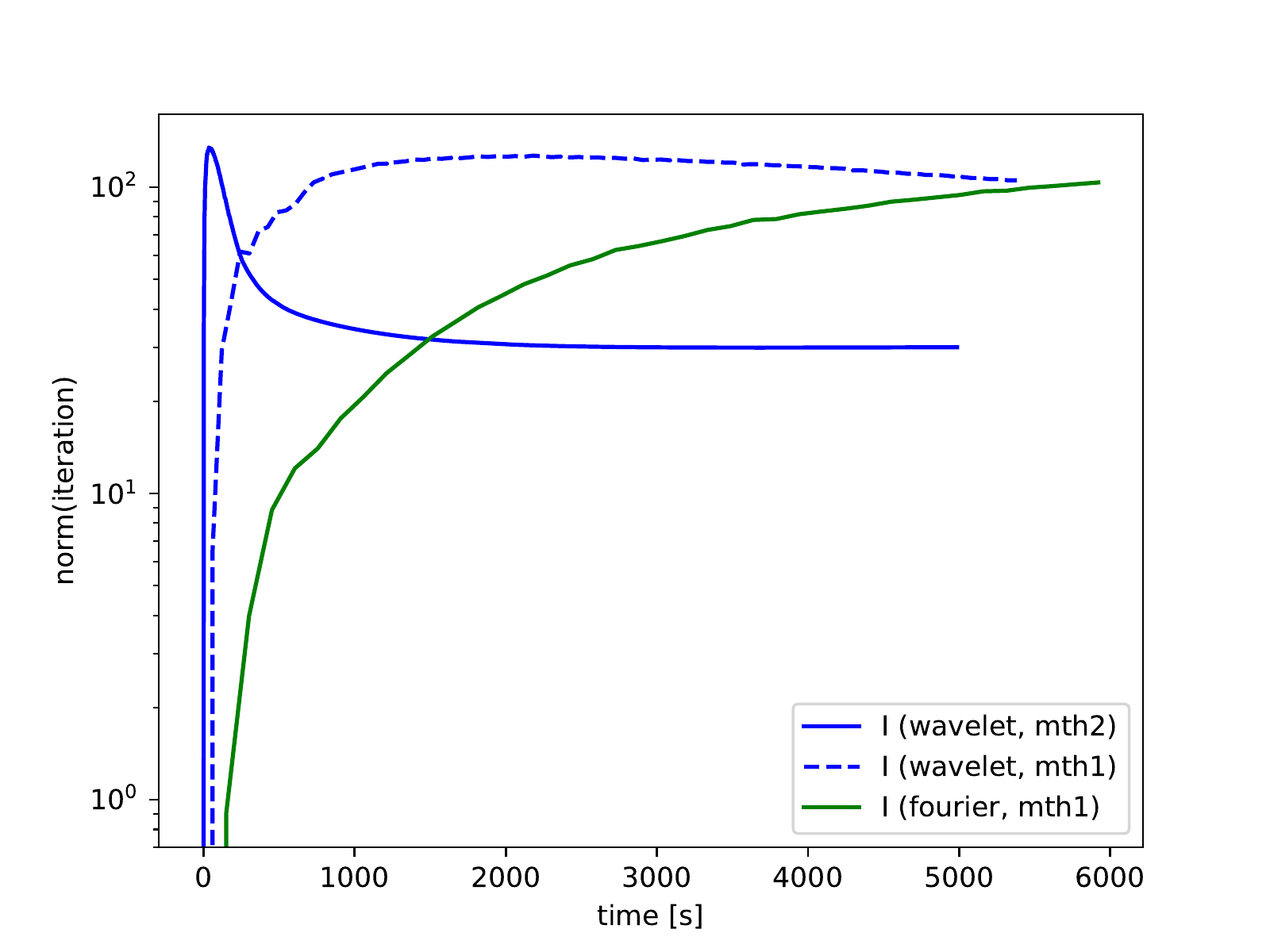}
\vspace*{2cm}
\caption{Plot of the Cameron-Martin-norm over time for all three relevant methods. This explains the abundance of detail in figure \ref{fig:wavelet_proposal} versus \ref{fig:wavelet_MAP}}
\label{fig:norms}
\end{minipage}
\end{figure}

\paragraph{Computation of sparse MAP estimators with a proper Besov prior}
The preceding simulations were presented with a choice of Gaussian wavelet prior, the main point of this paper being fast computation of the MAP estimator, not the sparsity of the reconstruction. Of course, the same results also hold for (non-Gaussian) Besov priors. See figure \ref{fig:sparse} for the MAP estimator of the inverse problem as described in the preceding section, but with a sparse Besov prior $B_1^{0.5}$ instead.

\begin{figure}[hbtp]
\centering
\begin{minipage}[t]{0.45\linewidth}
\includegraphics[width=\textwidth]{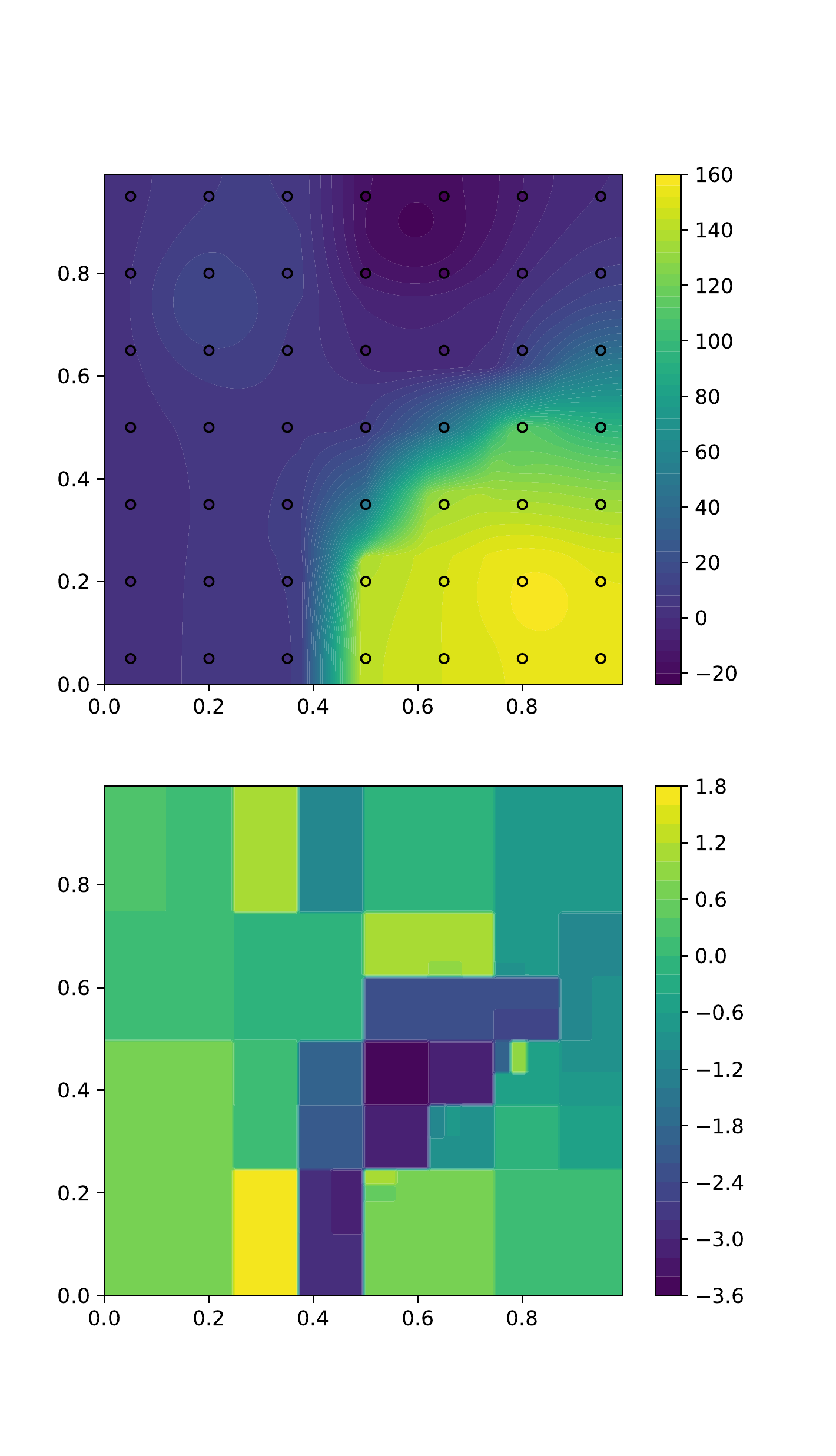}
\caption{MAP estimator (with method 2) for a Besov prior. }
\label{fig:sparse}
\end{minipage}
\hspace*{0.05\linewidth}
\begin{minipage}[t]{0.45\linewidth}
\includegraphics[width=\textwidth]{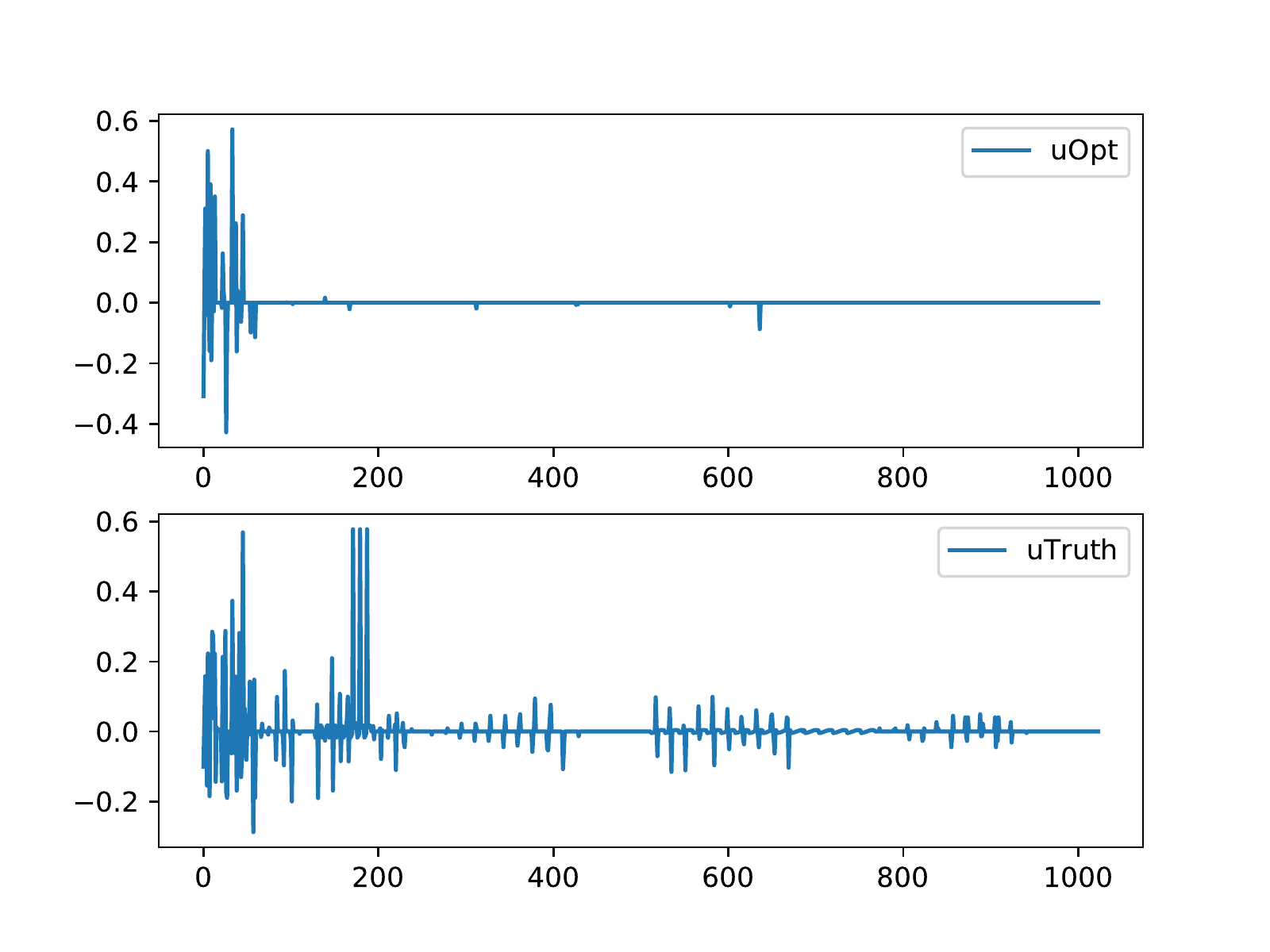}
\caption{Coefficients of the MAP estimator and ground truth. Note the much higher sparsity of the MAP estimator.}
\label{fig:coeffs_sparse}
\end{minipage}
\end{figure}

\section{Conclusion}
We have demonstrated how the use of wavelet-based priors facilitates much faster calculation of the misfit gradient needed for employment of more advanced optimization techniques like quasi-newton methods. This yields a major speed-up in the calculation of the MAP estimator of both the Gaussian wavelet prior $B_2^s$ but also of a proper Besov prior $B_1^s$ and can be used equivalently for sampling from the resulting posteriors if gradient-informed sampling algorithms like HMC \cite{girolami2011riemann}, MALA \cite{roberts1996exponential} or DILI \cite{cui2016dimension} are utilized.
}

\begin{appendix}
\section{Haar Wavelets}
The first wavelet was the Haar wavelet, conceived by Alfred Haar \cite{haar1910theorie} as a method of generating an unconditional basis for $L^2([0,1])$ (the Haar wavelet in fact builds an unconditional basis also for general $L^p$ spaces). Wavelets really took off in the 70s and 80s with notably Ingrid Daubechies' description of compactly supported continuous wavelets \cite{daubechies1988orthonormal} (the Haar wavelet is discontinous) and St\'ephane Mallat's general framework of multiresolution analysis \cite{mallat1989multiresolution}. A classical introduction are Daubechies' marvelous ``ten lectures"  \cite{daubechies1992ten}, a very instructive and readable account is by Blatter \cite{blatter2003wavelets}. In the following, we will only use Haar's original wavelets. The discontinuity (which can be a disadvantage for some applications e.g. in image processing) is actually a favorable property in our case as we want to model sharp interfaces in subsurface topology and thus we do not need or want more sophisticated wavelets (at least for the purpose of this study).
\subsection{Haar wavelet expansions of functions on \texorpdfstring{$[0,1]^d$}{[0,1]**d}}
We will work with dimension one and two; the step to higher dimensions is straightforward (although computationally challenging!).
\subsubsection{Dimension one} Define the Haar scale function and Haar mother wavelet in 1d by
\[ \phi(x) = \1{[0,1)}(x),\quad \psi(x) = \1{\left[0,\tfrac{1}{2}\right)}(x) - \1{\left[\tfrac{1}{2}, 1\right)}(x),\]
where $\1{A}(x) = 1$ if $x\in A$ and $\1{A}(x) = 0$ else.
The scaled wavelets are
\[ \psi_{j,k}(x) = 2^{j/2}\cdot \psi(2^{j}x-k).\]
We know that we can write any $f\in L^2[0,1]$ as 
\begin{align}
f(x) &= w_0\cdot \phi(x) + \sum_{j=0}^\infty \sum_{k=0}^{2^j-1}w_{j,k}\cdot \psi_{j,k}(x)\label{eq:wavelet1d_1} \\
&=w_0\cdot \phi(x) + \sum_{l=1}^\infty c_l \psi_l(x) \label{eq:wavelet1d_2} \end{align} 
where the second line is a re-indexing of $(j,k)\mapsto l(j,k) = 2^j+k$ and $\psi_l = \psi_{j,k}$ for the appropriate change of index (see figure \ref{fig:index1}). The justification for this to work is multiresolutional analysis (see \cite{daubechies1992ten}), but the main point is that wavelets give a way to give a \textit{local} way of spanning functions (as compared to, say, trigonometric polynomials -- the exact opposite of local) which is also \textit{robust with respect to fine discretization} (as opposed to a naive spanning of a function by evaluating it on a grid). We show how the wavelet expansion is computed for a concrete function:
\begin{figure}\label{fig:index1}
\begin{center}
\begin{tabular}{lll}
$l$ & $j$ & $k$\\
\hline
$1$ & $0$&$0$ \\ 
$2,3$ & $1$&$0,1$ \\ 
$\vdots$ & $\vdots$ &\\
$2^i,\ldots,2^{i+1}-1$ & $i$&$0,\ldots,2^i-1$\\
\end{tabular} 
\caption{Index table for dimension one}
\end{center}
\end{figure}
Let $f$ be given by its values on a uniform grid of $[0,1)$. We always assume that the number of grid points is a power of $2$, i.e. $\vec f = (f(0), f(2^{-N}),\ldots f(1-2^{-N}))^T$. Now we define $a_N^{(k)} \defeq f(k\cdot 2^{-N})$. Then iteratively we compute \begin{align*}
a_{j}^{(k)} &\defeq \frac{a_{j+1}^{(2k)}+a_{j+1}^{(2k+1)}}{2}\\
d_{j}^{(k)} &\defeq \frac{a_{j+1}^{(2k)}-a_{j+1}^{(2k+1)}}{2}
\end{align*} 
for $j=0,\ldots, N-1$ and $k = 0,\ldots, 2^n-1$, i.e. every step $j\mapsto j+1$ halves the size of the vectors $a_n$ and $d_n$ (by coarsening the resolution by a factor of two). Note that we can forget now all $a_j$ for $j>0$, as $a_{j+1}^{(2k)} = a_n^{(j)} + d_n^{(k)}$ and $a_{j+1}^{(2k+1)} = a_n^{(j)} - d_n^{(k)}$. Then the wavelet expansion is given by the form
\[ w_0\cdot \phi(x) + \sum_{j=0}^{N-1} \sum_{k=0}^{2^j-1}w_{j,k}\cdot \psi_{j,k}(x)\]
 with $w_0 = a_0$ and  $w_{j,k} = d_j^{(k)}$. Note that the sum is finite as having started with values of $f$ on a grid puts a lid on the maximal resolution.    
\begin{figure}
\begin{center}
\begin{tikzpicture}
\node (aN)[rectangle, draw=gray] {$a_N\,\,~$};
\node (aN-1) [right = of aN, rectangle, draw=gray] {$a_{N-1}$};
\node (dN-1) [below = of aN-1 , rectangle, draw=gray, fill=black!20] {$d_{N-1}$};
\node (aN-2) [right = of aN-1, rectangle, draw=gray] {$a_{N-2}$};
\node (dN-2) [below = of aN-2 , rectangle, draw=gray, fill=black!20] {$d_{N-2}$};
\node (dots) [right = of aN-2] {$\cdots$};
\node (a0) [right = of dots, rectangle, draw=gray, fill=black!20] {$a_{N-2}$};
\node (d0) [below = of a0 , rectangle, draw=gray, fill=black!20] {$d_{N-2}$};
\draw[->] (aN.east)-- (aN-1.west);
\draw[->] (aN-1.east)-- (aN-2.west);
\draw[->] (aN.east)-- (dN-1.west);
\draw[->] (aN-1.east)-- (dN-2.west);
\draw[->] (aN-2.east)-- (dots.west);
\draw[->] (dots.east)-- (a0.west);
\draw[->] (dots.east)-- (d0.west);
\end{tikzpicture}
\caption{Expansion calculation. Filled nodes represent quantities directly used for the expansion.}
\end{center}
\end{figure}
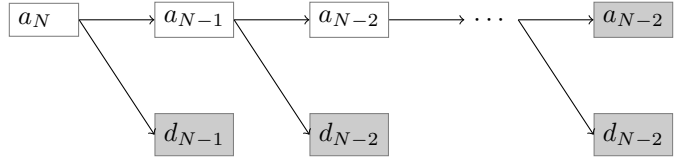

\begin{figure}[hbtp]
\centering
\includegraphics[width=0.65\textwidth]{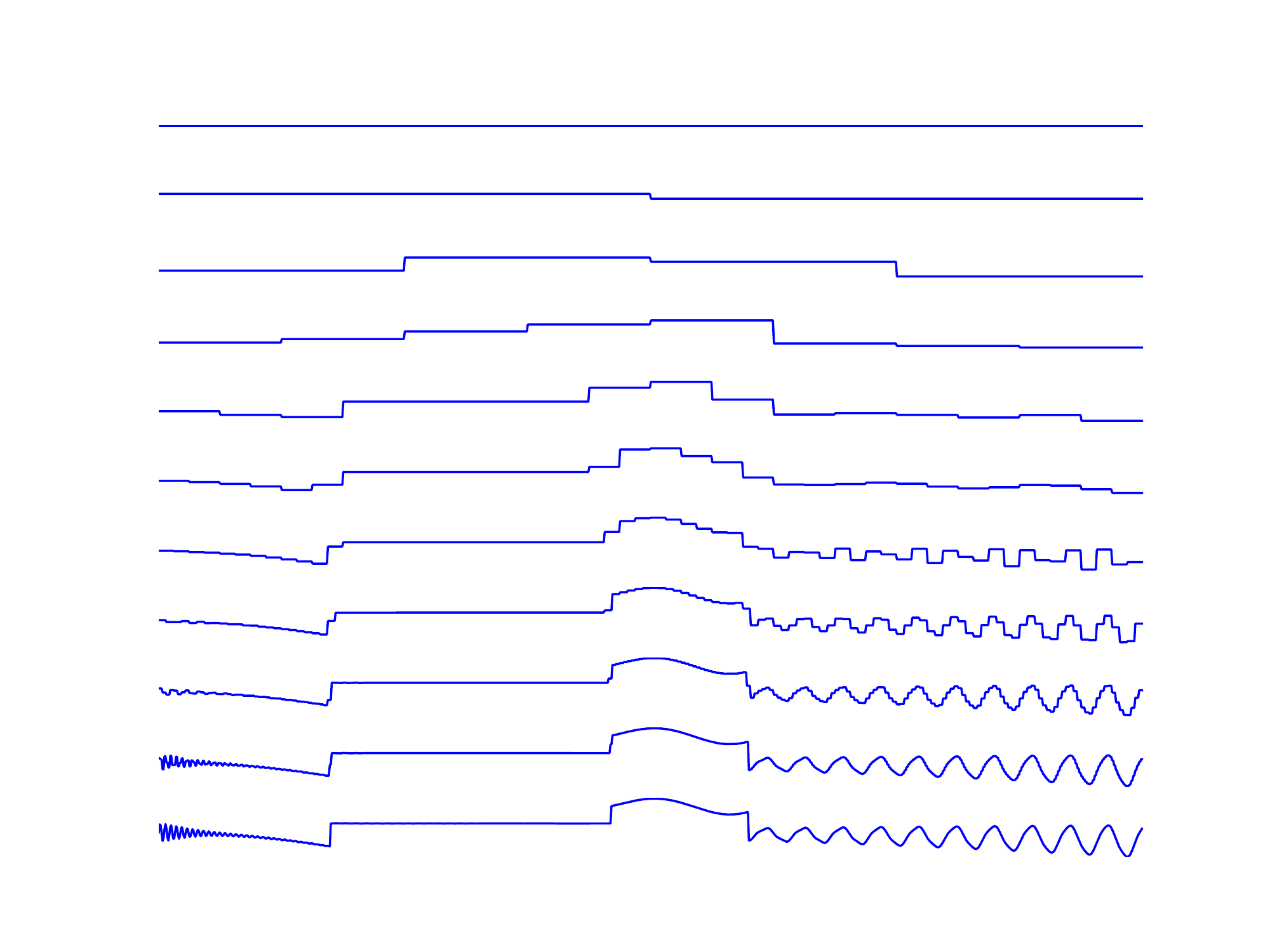}\includegraphics[width=0.35\textwidth]{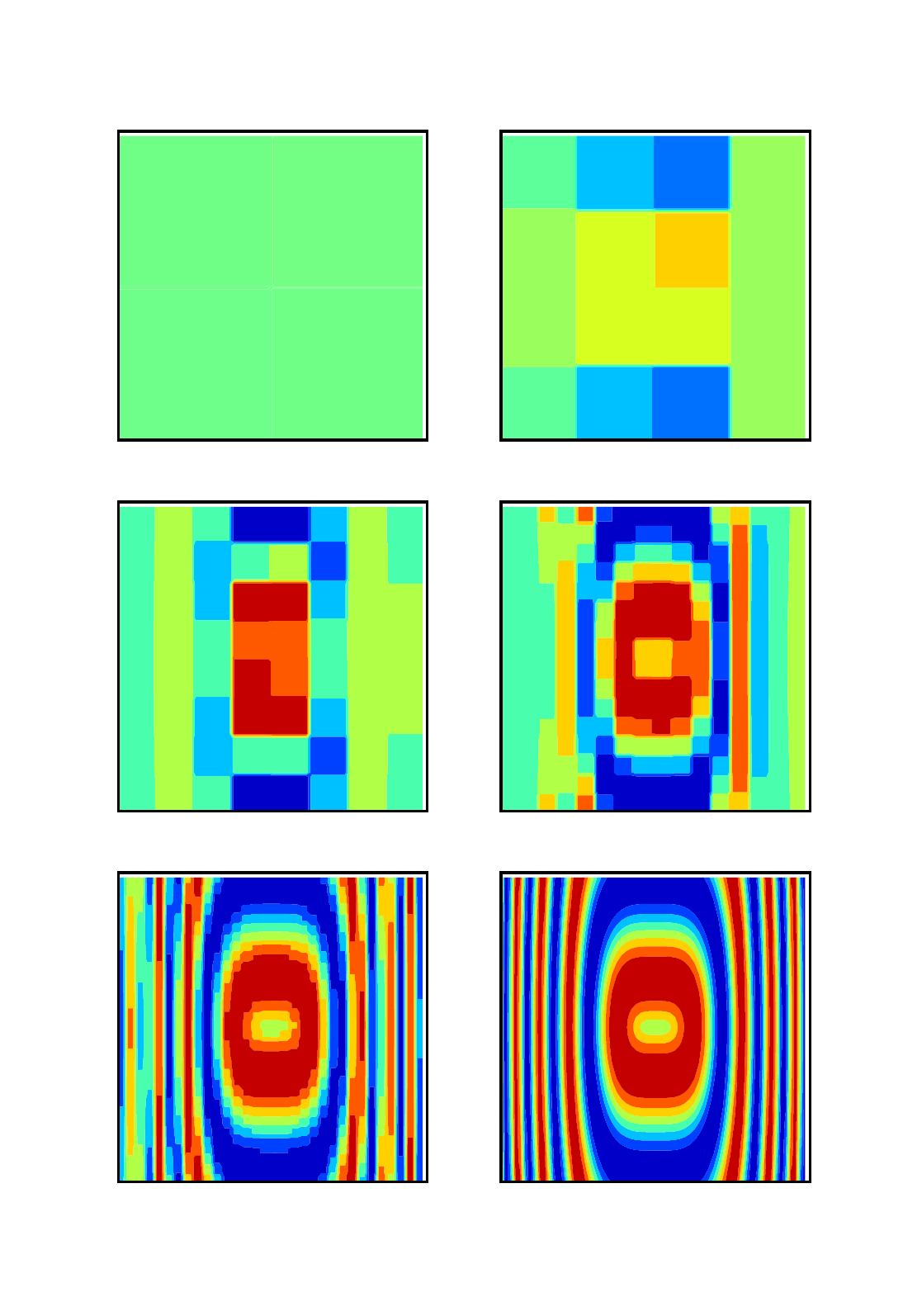}
\caption{Left: Reconstruction of a one-dimensional function from its wavelet decomposition: $\ell$-th row shows $a_0\cdot \phi(x) + \sum_{j=0}^{\ell-1}\sum_{k=0}^{2^j-1}d_j^{(k)}\cdot \psi_{j,k}(x)$. Right: Reconstruction of a two-dimensional function from its wavelet decomposition.}
\end{figure}

\subsubsection{Dimension two} 
Here we need to define the scale function and three mother wavelets.
\begin{align*}
\phi(x,y) &= \phi(x)\cdot\phi(y)\\
&=\1{[0,1)^2}(x,y)\\[1em]
\psi^{(0)}(x,y) &= \phi(x)\cdot\psi(y)\\
&=\1{[0,1)\times [0,1/2)}(x,y) - \1{[0,1)\times [1/2,1)}(x,y)\\[1em]
\psi^{(1)}(x,y) &= \psi(x)\cdot\phi(y)\\
&=\1{[0,1/2)\times [0,1)}(x,y) - \1{[1/2,1)\times [0,1)}(x,y)\\[1em]
\psi^{(2)}(x,y) &= \psi(x)\cdot\psi(y)\\
&=\1{[0,1/2)\times [0,1/2)}(x,y) - \1{[0,1/2)\times [1/2, 1)}(x,y)  \\
&- \1{[1/2, 1)\times [0,1/2)}(x,y)  + \1{[1/2, 1)\times [1/2, 1)}(x,y) 
\end{align*}
and scaled wavelets 
\[\psi_{j,k,n}^{(m)}(x,y) = 2^j\cdot \psi^{(m)}(2^jx-k, 2^jy-n).\]
With this we can expand a function defined on $[0,1]^2$ by
\begin{align*}
f(x,y) &= w_0 \cdot \phi(x,y) + \sum_{j=0}^\infty \sum_{m=0}^2\sum_{k=0}^{2^j-1}\sum_{n=0}^{2^j-1} w_{j,k,n}^{(m)}\cdot \psi_{j,k,n}^{(m)}(x,y)\\
&= w_0 \cdot \phi(x,y) + \sum_{l=1}^\infty c_l\psi_l(x,y)\end{align*} 
where as in 1d we re-index $(j, m, k, n)\mapsto l(j,m,k,n) = 4^j + m\cdot 4^j + k\cdot 2^j + l$ or as in the table in figure \ref{fig:index2}. 
\begin{figure}\label{fig:index2}
\begin{center}
\begin{tabular}{lllll}
$l$ & $j$ & $m$ & $k$ & $n$\\
\hline
\hline
$1$ & $0$ & $0$ &  $0$ &  $0$ \\
$2$ & $0$ & $1$ & $0$&  $0$ \\
$3$ & $0$ & $2$ & $0$&  $0$ \\ 
\hline
$4,\ldots 7$ & $1$ & $0$ & $0,1$& $0,1$\\
$8,\ldots 11$ & $1$ & $1$ & $0,1$ & $0,1$\\
$12,\ldots 15$ & $1$ & $2$ & $0,1$ & $0,1$\\
\hline
$16,\ldots 31$ & $2$ & $0$ & $0,\ldots 3$ & $0,\ldots 3$ \\
$32,\ldots 47$ & $2$ & $1$ & $0,\ldots 3$& $0,\ldots 3$\\
$48,\ldots 63$ & $2$ & $2$ & $0,\ldots 3$& $0,\ldots 3$\\
\hline
$4^i,\ldots 2\cdot 4^i-1$ & $i$ & $0$ & $0,\ldots 2^i-1$& $0,\ldots 2^i-1$\\
$2\cdot 4^i,\ldots 3\cdot 4^i-1$ & $i$ & $1$ & $0,\ldots 2^i-1$& $0,\ldots 2^i-1$\\
$3\cdot 4^i,\ldots 4^{i+1}-1$ & $i$ & $2$ & $0,\ldots 2^i-1$& $0,\ldots 2^i-1$\\
\end{tabular} 
\caption{Index table for dimension two}
\end{center}
\end{figure}
Calculation of the wavelet expansion in two dimension is done as follows: 

Given a function by its values on a square, power-of-$4$ grid $\{0, 2^{-N},\ldots 1-2^{-N}\}^2$ we define $a_N^{(k,n)} = f(k\cdot 2^{-N}, n\cdot 2^{-N})$ and
\begin{align*}
a_j^{(k,n)} &\defeq \frac{a_{j+1}^{(2k,2n)}+a_{j+1}^{(2k+1,2n)}+a_{j+1}^{(2k,2n+1)}+a_{j+1}^{(2k+1,2n+1)}}{4}\\
d_j^{0, (k,n)} &\defeq \frac{a_{j+1}^{(2k,2n)}+a_{j+1}^{(2k+1,2n)}-a_{j+1}^{(2k,2n+1)}-a_{j+1}^{(2k+1,2n+1)}}{4}\\
d_j^{1, (k,n)} &\defeq \frac{a_{j+1}^{(2k,2n)}-a_{j+1}^{(2k+1,2n)}+a_{j+1}^{(2k,2n+1)}-a_{j+1}^{(2k+1,2n+1)}}{4}\\
d_j^{2, (k,n)} &\defeq \frac{a_{j+1}^{(2k,2n)}-a_{j+1}^{(2k+1,2n)}-a_{j+1}^{(2k,2n+1)}+a_{j+1}^{(2k+1,2n+1)}}{4}
\end{align*}
and analogously to one dimension, we can then span $f$ by 
\[a_0 \cdot \phi(x,y) + \sum_{j=0}^{N-1} \sum_{m=0}^2\sum_{k=0}^{2^j-1}\sum_{n=0}^{2^j-1} d_j^{m,(k,n)}\cdot \psi_{j,k,n}^{(m)}(x,y)\]
where the meaning of the indices is as follows:
\begin{itemize}
\item $j$: Scale
\item $m$: Orientation ($0$=Horizontal, $1$=Vertical, $2$=Diagonal)
\item $k$: Shift in horizontal direction
\item $n$: Shift in vertical direction
\end{itemize}

\subsubsection{Arbitrary dimension}
Just for reference we give the straight-forward extension to arbitrary dimension $d\in\N$: Given a function $f :[0,1]^d\to \R$, we define $\Lambda = \{0,\ldots,2^j-1\}^d$ and span
\[f(z) = w_0 \cdot \phi(z) + \sum_{j=0}^\infty \sum_{m=0}^{2^d-1}\sum_{k\in \Lambda} w_{j,k}^{(m)}\cdot \psi_{j,k}^{(m)}(z)\]

\section{Besov spaces}

We can characterize (see \cite{triebel2008function,daubechies1992ten,meyer1995wavelets,bui2015scalable,dashti2011besov,kolehmainen2012sparsity,lassas2009discretization}) elements of Besov space $B_{pp}^s$ by the following: Let $f:[0,1]^d\to\R$. Assume that the wavelet basis used is regular enough. Then for $p>0$ and $s \in\R$,
\begin{equation}\label{eq:besovdimindep}
 f\in B_{pp}^s([0,1]^d) \Leftrightarrow \left( \sum_{l=1}^\infty l^{\frac{ps}{d}+\frac{p}{2}-1}\cdot |c_l|^p\right)^\frac{1}{p} < \infty. \end{equation}
In one and two dimensions this reduces to the following:

 \begin{align}
 f&\in B_{pp}^s([0,1]), \, d=1\nonumber \\
  &\Leftrightarrow \|f\|_{B_{pp}^s} \defeq \left(|w_0|^p + \sum_{j=0}^\infty 2^{jp(s+\frac{1}{2}-\frac{1}{p})}\cdot \sum_{k=0}^{2^j-1} |w_{j,k}|^p \right)^\frac{1}{p} < \infty \label{eq:besovdim1} \\
 f&\in B_{pp}^s([0,1]^2), \, d=2 \nonumber\\
  &\Leftrightarrow \|f\|_{B_{pp}^s} \defeq \left(|w_0|^p + \sum_{j=0}^\infty 2^{jp(s+1-\frac{2}{p})}\cdot \sum_{m=0}^2\sum_{k=0}^{2^j-1}\sum_{n=0}^{2^j-1} |w_{j,k,n}^{(m)}|^p \right)^\frac{1}{p} < \infty.\label{eq:besovdim2} 
 \end{align}

Note that the formula in \eqref{eq:besovdimindep} does not yield the same numerical value as the appropriate dimension-dependent formula for $\|f\|_{B_{pp}^s}$ but that the expressions are merely equivalent in the sense of equivalent norms. 

If we set $p=2$, we recover the Sobolev spaces, i.e. $B_{2,2}^s = H^s$.
\begin{rem}
Note that Besov spaces usually have an additional parameter, i.e. $B_{p,q}^s$ although we only use the case $p=q$ for our choice of priors. For completeness, (both for general $p>0,q>0,s\in\R$ and dimension $d$) the $B_{p,q}^s$ Besov norm for functions defined on $[0,1]^d$ is defined by
\[\|f\|_{B_{pq}^s([0,1]^d)}\defeq \left(|w_0|^p + \sum_{j=0}^\infty 2^{jq(s+\frac{d}{2}-\frac{d}{p})}\left(\sum_{m=0}^{2^d-1}\sum_{k\in\Lambda}|w_{j,k}^{(m)}|^p\right)^\frac{q}{p} \right)^\frac{1}{q}\]
\end{rem}
\end{appendix}

\section*{Acknowledgments} P.W. is thankful for a fruitful discussion with Youssef Marzouk, a very helpful email from Donald Estep and in particular for the guidance of Claudia Schillings which ultimately led to the idea of this paper.
\bibliographystyle{siam}
\bibliography{lit}

\end{document}